\documentclass[twoside, twocolumn]{article}
\usepackage[a4paper, margin=2.236068cm]{geometry}
\usepackage[LGR, T1]{fontenc}
\usepackage{amsmath}
\usepackage{amssymb}
\usepackage{amsfonts}
\usepackage{mathtools}
\usepackage{booktabs}
\usepackage{nicefrac}
\usepackage{graphicx}
\usepackage{tabularx}
\usepackage{float}
\usepackage{chngcntr}
\usepackage{fancyhdr}
\usepackage{titlesec}
\usepackage{titletoc}
\usepackage{emptypage}
\usepackage{tocloft}
\usepackage{bm}
\usepackage{indentfirst}
\usepackage{chngcntr}
\usepackage[justification=centering]{caption}
\usepackage[style=numeric-comp, sorting=none]{biblatex}
\usepackage{hyperref}
\usepackage{breakurl} 
\usepackage[nameinlink]{cleveref}
\usepackage{titling}
\usepackage{ifoddpage}
\usepackage{graphicx}
\usepackage{stfloats}
\usepackage{adjustbox}
\usepackage{tabulary}
\usepackage{soul}

\newcommand{\mytitlepartone}{Universal Definitions of the Roman Factorial}
\newcommand{\mytitleparttwo}{Introduction to Foundational Functions and the Generalization Process}
\newcommand{\spm}{\vcenter{\hbox{\scalebox{0.868034}{$\pm$}}}}
\newcommand{\bp}[1]{\big( #1 \big)}
\newcommand{\bb}[1]{\big[ #1 \big]}
\newcommand{\cf}[1]{\lceil #1 \rceil}
\newcommand{\ff}[1]{\lfloor #1 \rfloor}
\newcommand{\rf}[1]{\lfloor #1 \rceil !}

\newcommand{\lrrf}[1]{\left\lfloor #1 \right\rceil}
\newcommand{\mf}[2]{#1!_{(#2)}}
\newcommand{\of}[2]{#1!_{#2}}

\newcommand{\mbR}{\mathbb{R}}
\newcommand{\mbZ}{\mathbb{Z}}
\newcommand{\mbN}{\mathbb{N}}
\newcommand{\detnref}[2]{
    \begin{flushleft}
        \large\textbf{#1}
        \addcontentsline{toc}{subsection}{#1}
    \end{flushleft}
    \markright{\MakeUppercase{#1}}

    \par #2
}
\let\|=\enspace
\newcommand{\M}{\quad}

\hypersetup{
  colorlinks=true,
  linkcolor=black,
  citecolor=black,
  filecolor=black,
  urlcolor=black,
  breaklinks=true,
  bookmarksopen=true,
  bookmarksnumbered=true,
  pdftitle={\mytitlepartone \text{: } \mytitleparttwo},
  pdfauthor={Leonidas Liponis},
  pdfsubject={Roman factorial generalization},
  pdfkeywords={factorial, roman factorial, generalization},
  pdfstartview={FitH},
  pdfpagemode={UseOutlines},
  pdfdisplaydoctitle=true
}
\makeatletter
\DeclareRobustCommand{\bfseries}{
  \not@math@alphabet\bfseries\mathbf
  \fontseries\bfdefault\selectfont
  \boldmath
}
\makeatother
\crefname{section}{Section}{Sections}
\crefname{subsection}{Subsection}{Subsections}
\crefname{subsubsection}{Subsubsection}{Subsubsections}
\DeclareCaptionLabelFormat{nofigure}{#2}
\DeclareCaptionLabelFormat{notable}{#2}
\crefname{table}{}{}
\crefname{figure}{}{}
\crefname{equation}{Eq.}{Eqs.}

\renewcommand{\thetable}{\text{Tbl.} \arabic{section}.\arabic{subsection}.\arabic{table}}

\makeatletter
\@addtoreset{table}{subsubsection}
\@addtoreset{figure}{subsubsection}
\@addtoreset{equation}{subsubsection}
\makeatother

\fancypagestyle{plain}{
  \fancyhf{}
  \fancyhf[EHL]{\thepage}
  
}
\let\oldsection\section
\renewcommand{\section}[1]{\vspace*{\fill}\pagebreak\oldsection{#1}}
\newcommand{\firstsection}[1]{\oldsection{#1}}

\title{\mytitlepartone \text{:} \mytitleparttwo}
\author{Leonidas Liponis}
\date{March 14, 2024}

\begin{document}

\twocolumn[
    \begin{@twocolumnfalse}
        \maketitle
        \begin{abstract}
            This paper introduces a new method for redefining the Roman factorial using universally applicable functions that are not expressed in closed form. We present a set of \textit{foundational functions}, similar to Boolean operations, to simplify the factorial expression. Through a systematic process of generalization, termed \textit{generalization process}, we aim to use these \textit{foundational functions} to create recursive and non-recursive, global definitions of the Roman factorial.
        \end{abstract}
        \vspace{1.618034em}
    \end{@twocolumnfalse}
]

\firstsection{Preface}\label{s:preface}

\textit{Universal Definitions of the Roman Factorial} is a study that expands the concept of the factorial to various domains and develops mathematical expressions to describe it.

The study is too large and extensive to be contained in a single paper, so it is divided in parts. \href{https://arxiv.org/abs/2403.09581}{Part 1} is considered an introduction to its follow-ups, but it stands alone as a analysis into the Roman factorial.

\subsection{Contents}\label{ss:contents}

The contents of this study are divided into 5 parts, of which this paper is the first one. Each part contains several sections, and they are further subdivided into subsections. For example, you are currently reading \href{https://arxiv.org/abs/2403.09581}{Part 1}, \Cref{s:preface}, \Cref{ss:contents}.

\href{https://arxiv.org/abs/2403.09581}{Part 1} contains four main sections. \Cref{s:basic ff,s:advanced ff} are about \textit{foundational functions} (\textit{F.F.}) while \Cref{s:rf rec rl gen,s:rf prod gen} are about redefining the Roman factorial recursively and non-recursively\footnote{The term "recursively", as well as the concept of "\textit{foundational functions}", will be clarified shortly. Any unfamiliar mathematical terminology introduced in this section or later will be explained in time.}. This paper focuses on introducing the \textit{F.F.} step by step, in a way that builds intuition and insight.

These \textit{F.F.} are used for expressing the Roman factorial definition in a single formula (instead of a 2-part piece-wise definition) and in constructing a universal, non-recursive formula across all integers. The procedure of incorporating \textit{F.F.} into expressions of the Roman factorial is referred to as the \textit{generalization process}, and will be followed in \Cref{s:rf rec rl gen,s:rf prod gen} as well as future papers.

Later parts of this study\footnote{The next part of this study will introduce a notation system for the different types of factorial variations. However, since this paper focuses exclusively on the Roman factorial, there's no need to distinguish between similar concepts here. Therefore, what might later be called $\mf{n}{1}$ or $\of{n}{1}$ will simply be referred to as $\rf{n}$ in this paper.} include applications of this \textit{generalization process} to the double factorial, triple factorial and all higher order factorials. There is a part about generalizing the factorial to non-integers, and another about non-integer orders. The last is about complex factorials and a new formulation of the \textit{F.F.} presented in this paper to accommodate for the new domain.

Concepts like those mentioned in the paragraph above will not be explained and analyzed here, so there is no need to learn about them now. In fact, \ul{no prior knowledge is required} about any advanced mathematical concept, as this paper was written with accessibility in mind. Anyone with high-school education should be able to read this study all the way through, which is its main goal in the first place.

For those interested in diving deeper into specific mathematical concepts, an addendum is provided in \Cref{s:addendum}. It has basic information about a variety of things and is intended to fill any gaps not explicitly elaborated upon in the main text. Lastly, \Cref{s:references} lists all references used in the formulation of this paper along with other details and information.

\subsection{Useful information}\label{ss:useful information}

If you are accessing this document online, it's worth noting that all references to equations, figures, tables, and sections are clickable. This feature allows for instant navigation to the referenced content without the need to remember anything.

Equations or tables may be repeated throughout this paper, but their original numbering is preserved. For instance, \ref{eq:rf pw} appears many times in different sections, but each time it has the same tag. It is also clickable as a hyperlink, referencing its first appearance in the text.

This paper uses terminology that may be unfamiliar at first. For this reason, a list of key words and phrases are provided here:

\begin{itemize}
    \item \textbf{Factorial}: A mathematical function for natural numbers, expressed as the product of positive integers.
    \item \textbf{Roman factorial}: An extension of the factorial into negative integers.
    \item \textbf{Foundational functions}: A set of Boolean-like functions introduced in this paper, each designed with straightforward outputs.
    \item \textbf{Recursive relationship}: A mathematical definition that expresses a sequence or function in terms of its own previous values through a specified formula or rule. Also: \textit{recurrence relation}
    \item \textbf{Non-recursive relationship}: A mathematical definition that directly expresses a sequence or function without referring to its own previous values. Also: \textit{direct, iterative}
    \item \textbf{Piece-wise definition}: A method of specifying a function or expression by dividing its domain into distinct intervals, each associated with a specific formula or rule. Also: \textit{closed form}
    \item \textbf{Universal definition}: A single mathematical expression of a function, applied uniformly across all input ranges. Also: \textit{unified, global}
    \item \textbf{Universal definitions of the Roman factorial}: Expressions of the Roman factorial that are described by a single formula across all integers.
    \item \textbf{Generalization process}: A way of inserting \textit{F.F.} into piece-wise definitions of a function, aiming to unite the cases into a single expression.
\end{itemize}

The abbreviations used throughout this paper are referenced here:
\begin{itemize}
    \item \textit{F.F.}: \textit{Foundational functions}, a series of simple Boolean-like functions that are built upon each other and usually have binary outputs (0 or 1).
    \item \textit{Eq.}: An equation or a relationship.
    \item \textit{Tbl.}: A mathematical table that will usually be a list of inputs and outputs of various functions.
    \item \textit{Fig.}: A figure or a diagram, that here will often depict the behavior of a function in a domain close to 0.
\end{itemize}

That concludes all introductory information. In the following section we'll introduce the factorial and an expansion to negative integers. If you encounter any unfamiliar terms or abbreviations, you can always refer back to these pages for clarification.

\section{Introduction to the factorial}\label{s:tf intro}
\subsection{Factorial}\label{ss:tf}

The factorial $(n!)$ is a mathematical function for natural numbers. It appears mainly in combinatorics, although it finds applications across various branches of mathematics. Specifically, $n!$ represents the product of all positive integers less than or equal to $n$. For example:
\begin{gather*}
    3! = 3 \cdot 2 \cdot 1 = 6 \,, \\[0.618034em]
    4! = 4 \cdot 3 \cdot 2 \cdot 1 = 24 \,, \\[0.618034em]
    5! = 5 \cdot 4 \cdot 3 \cdot 2 \cdot 1 = 120 \,, \\ \vdots \\[0.618034em]
    n! = n(n - 1)(n - 2) \cdots \,3 \cdot 2 \cdot 1 \,.
\end{gather*}

The factorial can be represented in two forms: as a rising/ascending or as a falling/descending product\footnote{Not to be confused with falling factorials or rising factorials, which are different concepts despite their similar names. More about these in \hyperref[ss:add fall rise fact]{Addendum~\ref*{ss:add fall rise fact}}.}. In the examples provided above, it's depicted as a falling product because each factor (or multiplicand) decreases by 1 with each step. A rising product looks like this:
\begin{equation*}
    3! = 1 \cdot 2 \cdot 3 \,.
\end{equation*}

Whether we conceptualize the factorial as an ascending or descending product, the final outcome remains unchanged. The only distinction lies in the arrangement of the factors within the product.

We can express the factorial as a $\prod$-product\footnote{More info about $\prod$-products in \hyperref[ss:add product]{Addendum~\ref*{ss:add product}}.}. Expressed as a rising product, it is
\begin{equation}\label{eq:tf prod rising}
    n! = \prod_{k = 1}^{n} k \,, \M n \in \mbZ^+ .
\end{equation}

If we consider the factorial as a descending product instead, we have
\begin{equation}\label{eq:tf prod falling}
    n! = \prod_{k = 0}^{n - 1} (n - k) \,, \M n \in \mbZ^+ .
\end{equation}

\ref{eq:tf prod rising} and \ref{eq:tf prod falling} are non-recursive definitions of the factorial. This means that values of $n!$ can be calculated without knowing the value of other factorials. For example, the value of $8!$ can be found directly using this definition: $8! = 40320$. Non-recursive definitions are also called direct.

Below, we list the first few factorials:
\begin{table}[H]
    \centering
    \begin{tabular}{c|*{8}{c}}
        $n$ & 0 & 1 & 2 & 3 & 4 & 5 & 6 & 7 \\
        \hline
        $n!$ & 1 & 1 & 2 & 6 & 24 & 120 & 720 & 5040 \\
    \end{tabular}
    \caption{Factorials}
    \label{tbl:f values}
\end{table}

We see that $0! = 1$. This is a convention for the empty product\footnote{More info about the empty product in \hyperref[ss:add product]{Addendum~\ref*{ss:add product}}.} and also because it extends various combinatoric identities. For more reasons why $0!$ should equal 1 see \hyperref[ss:add 0 factorial]{Addendum~\ref*{ss:add 0 factorial}}, but we will show a few of them shortly.

Let's now investigate the recursiveness of the factorial. In the following examples, we list a few factorials again. It appears that $n!$ equals the product of the number $n$ with the next smaller factorial:
\begin{gather*}
    3! = 3 \cdot 2! = 3 \cdot 2 \cdot 1 \,, \\[0.618034em]
    4! = 4 \cdot 3! = 4 \cdot 3 \cdot 2 \cdot 1 \,, \\[0.618034em]
    5! = 5 \cdot 4! = 5 \cdot 4 \cdot 3 \cdot 2 \cdot 1 \,, \\ \vdots \\[0.618034em]
    n! = n (n - 1)! \,.
\end{gather*}

We can include the starting value $0! = 1$ to define the factorial recursively, like so:
\begin{equation}\label{eq:tf rec}
    n! = n(n - 1)! \,, \M 0! = 1 \,, \M n \in \mbZ^+ .
\end{equation}

\ref{eq:tf rec} is called the recursive definition of the factorial, because in order to calculate $n!$ we need to know $(n - 1)!\,$, starting with $0! = 1$. Using this relationship we can find the next larger factorial for a given factorial. For example, if we know what $8!$ is, we can easily find $9!$ by evaluating for $n = 9$: $9! = 9 \cdot 8! = 362880$.

By defining $0!$ to be 1, we can expand \ref{eq:tf rec} to include the $n = 0$ case, thus rendering it concise.

If we substitute $n$ with $n + 1$ in \ref{eq:tf rec} and solve for $n!\,$, we get a similar equation\footnote{In this paper, the symbol $\mbZ^+$ refers to natural numbers or positive integers $(1, 2, 3, \dots)$ while $\mbZ^+_0$ means non-negative integers $(0, 1, 2, 3, \dots) \,$. More information about number sets in \hyperref[ss:add number sets]{Addendum~\ref*{ss:add number sets}}.}:
\begin{equation}\label{eq:tf rec prev}
     n! = \cfrac{(n + 1)!}{n + 1} \;, \M n \in \mbZ^+_0 \,.
\end{equation}

This relationship, although it does not define the factorial recursively, is useful for finding the value of a factorial when the next larger factorial is known. For example, let's find the factorial when $n = 0$:
\begin{gather*}
    3! = \cfrac{4 \cdot 3 \cdot 2 \cdot 1}{4} = \cfrac{4!}{4} \:, \\[0.618034em]
    2! = \cfrac{3 \cdot 2 \cdot 1}{3} = \cfrac{3!}{3} \:, \\[0.618034em]
    1! = \cfrac{2 \cdot 1}{2} = \cfrac{2!}{2} \:, \\[0.618034em]
    0! = \cfrac{1}{1} = 1 \,.
\end{gather*}

This demonstrates another reason why $0! = 1$.

\begin{figure}[H]
    \centering
    \includegraphics[width = 1\linewidth]{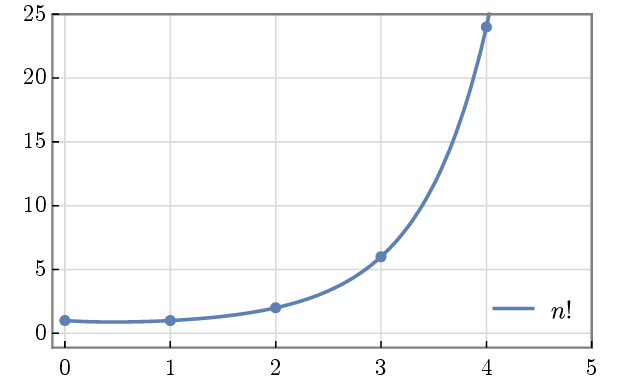}
    \caption{Factorial}
    \label{fig:factorial}
\end{figure}

The factorial is depicted below in \Cref{fig:factorial}. Each point corresponds to integer factorials, while the Gamma function is used to connect these discrete points smoothly. This continuation will be presented in later parts of this study: its role in \Cref{fig:factorial} is solely to be visually pleasing.

When trying to find factorials of negative numbers, for example $(-1)!\,$, we come across a problem. We initially notice that there is no intuitive definition of the factorial for negative integers, so we need a new definition entirely or at least an expansion of the recursive relationship as iterated in \ref{eq:tf rec prev}.

Let's approach this issue by attempting to find $(-1)!$ recursively. Using the recurrence relationship in \ref{eq:tf rec prev} for $n = -1$, we observe that we have to divide by 0:
\begin{equation*}
    (-1)! = \cfrac{(-1 + 1)!}{-1 + 1} = \cfrac{0!}{0} = \cfrac{1}{0} \;,
\end{equation*}
and that is undefined.

Also, if we use \ref{eq:tf rec} for $n = 0$, we see that
\begin{equation*}
    0! = 0 \cdot (-1)! \|\Rightarrow\|1 = 0 \cdot (-1)! \,,
\end{equation*}
which is impossible since there is no value of $(-1)!$ that satisfies this equation (no number multiplied by 0 gives 1).

Furthermore, if we extend this process, $(-2)!$ is undefined and so are all negative integer factorials. If we attempt to evaluate them recursively, we come across other undefined values:
\begin{align*}
    (-2)! &= \cfrac{(-2 + 1)!}{-2 + 1} = \cfrac{(-1)!}{-1} \:, \\[0.618034em]
    (-3)! &= \cfrac{(-3 + 1)!}{-3 + 1} = \cfrac{(-2)!}{-2} \:\cdots
\end{align*}

In the next subsection we will examine an extension of the factorial to negative integers. Note that the extension does not unveil any hidden behavior of the factorial function, as there are no traditional negative factorials. Rather, it is an expansion of the definition that has uses in some contexts.

\subsection{Roman factorial}\label{ss:rf}

The Roman factorial \cite{roman} \cite{romanfact} is an extension of the factorial function to negative integers\footnote{The symbols \cite{roman} and \cite{romanfact} are references to citations, listed in \Cref{s:references}.}. It is named after Steven Roman, who used it in his work about umbral calculus. Roman factorial can be expressed in the following closed form:
\begin{equation}\label{eq:rf pw}
    \rf{n} = 
    \begin{cases}
        \,n! & , \| n \in \mbZ^+_0 \\
        \cfrac{(-1)^{-n - 1}}{(-n - 1)!} & , \| n \in \mbZ^- ,
    \end{cases}
\end{equation}
in which the factorial can be defined recursively as
\begin{equation}\tag{\ref{eq:tf rec}}
    n! = n(n - 1)! \,, \M 0! = 1 \,, \M n \in \mbZ^+ .
\end{equation}

The variable $n$ can take only integer values, as denoted by the sets $\mbZ^+_0$ and $\mbZ^-$. The following table is a list of negative integer Roman factorials, as defined by \ref{eq:rf pw}:
\begin{table}[H]
    \centering
    \begin{tabular}{c|*{7}{c}}
        $n$ & -7 & -6 & -5 & -4 & -3 & -2 & -1 \\
        \hline
        $\rf{n}$ & $\nicefrac{1}{720}$ & -$\nicefrac{1}{120}$ & $\nicefrac{1}{24}$ & -$\nicefrac{1}{6}$ & $\nicefrac{1}{2}$ & -1 & 1
    \end{tabular}
    \caption{Negative Roman factorials}
    \label{tbl:rf values}
\end{table}

Negative Roman factorials exhibit a pattern akin to the reciprocals of positive factorials, with alternating signs. Notably, there is an offset: $\rf{-5}$ is not the reciprocal of $5!$, but rather of $4!$. However, this offset does not affect any properties of the factorial.

This definition obeys the recursive relationships of \ref{eq:tf rec} or \ref{eq:tf rec prev} in the entirety of their domain, which is now all integers. The latter of the relationships can now be written like this:
\begin{equation}\label{eq:rf rec prev}
     \lrrf{n}! = \cfrac{\lrrf{n + 1}!}{n + 1} \;, \M n \in \mbZ^- \setminus \{-1\}.
\end{equation}

This relationship is defined for any integer $n$ except -1, which is expressed mathematically by the set\footnote{More about number sets on \hyperref[ss:add number sets]{Addendum~\ref*{ss:add number sets}}.} $\mbZ^- \setminus \{-1\}$. \ref{eq:rf rec prev} does not define the Roman factorial recursively, it is just an identity valid for all integers except $-1$.

Alternatively, the Roman factorial can be defined recursively in two domains, like so:
\begin{equation}\label{eq:rf rec pw}
    \rf{n} = 
    \begin{cases}
        n \rf{n - 1} & , \| n \in \mbZ^+ \\[0.368034em]
        1 & , \| n = \{0 \,, -1\} \\[0.118034em]
        \cfrac{\rf{n + 1}}{n + 1} & , \| n \in \mbZ^- \setminus \{-1\} \,.
    \end{cases}
\end{equation}

Let's confirm \ref{eq:rf rec pw}'s applicability to negative integers. Using $\rf{-2}$ as a example, we can see that there is no issue in the calculation:
\begin{equation*}
     \rf{-2} = \cfrac{\rf{-2 + 1}}{-2 + 1} = \cfrac{\rf{-1}}{-1} = -1 \,.
\end{equation*}

The value of $(-2)!$ would normally be found using the recursive relationships of the normal factorial, if only $(-1)!$ was properly defined. $\rf{-1}$ has a well-defined value only in the context of the Roman factorial, which provides values for all negative integers while being recursively consistent. Additionally, the division by $n + 1$ is the reason for the alternating signs of negative integer factorials, since that divisor is negative when $n < -1$.

Certainly, there remains an anomaly in \ref{eq:rf rec prev} particularly at $n = -1$, yet the piece-wise Roman factorial offers a solution: $\rf{-1} = 1$. This specific value doesn't resolve the division by 0; indeed, no number can. Instead, the recursive relationship fails to provide a value for $(-1)!$ or $\rf{-1}$, so we address this limitation by asserting that the relationship is undefined at this anomaly (assigning $(-1)!$ the value of 1 arbitrarily).

In these last few paragraphs of the introduction, we will show an example of generalizing \ref{eq:rf rec pw} so that it becomes a single expression, universally defined across its domain. To begin with, consider the following Boolean-valued functions\footnote{More about Boolean functions in \hyperref[ss:add boolean algebra]{Addendum~\ref*{ss:add boolean algebra}}.}:
\begin{gather}
    f_1 \equiv f_1(n) = 
    \begin{cases}\label{eq:boolean f1 intro}
        +1 \,, & n \geq 0 \\
        -1 \,, & n < 0 \,,
    \end{cases}
    \\f_2 \equiv f_2(n) = 
    \begin{cases}\label{eq:boolean f2 intro}
        0 \,, & n \geq 0 \\
        +1 \,, & n < 0 \,.
    \end{cases}
\end{gather}

Using these, we can rewrite \ref{eq:rf rec pw} as follows:
\begin{equation}\label{eq:boolean rf preview}
    \begin{gathered}
        \rf{n} = \bp{n + f_2}^{f_1} \rf{n - f_1} \,, \| n \neq \{0, -1\} \,, \\[0.618034em]
        \rf{0} = \rf{-1} = 1 \,, \M n \in \mbZ \,.
    \end{gathered}
\end{equation}

The process of reaching this result will be explained thoroughly in \Cref{s:rf rec rl gen} and it serves as an example of the \textit{generalization process} used in this paper. This procedure seeks to rewrite various closed-form definitions as singular expressions by employing a set of \textit{foundational functions} which have simple and analytic forms.

\Cref{s:rf rec rl gen,s:rf prod gen} aim to establish two universal definitions for integer Roman factorials. One approach will result in a recursive definition and the other will result in a direct one, both applicable to all integers.

Lastly, the non-recursive definition of the Roman factorial will have the form
\begin{equation}\label{eq:rf gen example}
    \rf{n} = F(n) \,, \M n \in \mbZ \,,
\end{equation}
where $F(n)$ stands for a product-like function, similar to \ref{eq:tf prod rising} or \ref{eq:tf prod falling} for positive integers. The recursive expression for the Roman factorial will be based on \ref{eq:rf rec pw} or \ref{eq:rf pw} but it will be a singular expression, universally defined across all integers.

To establish these definitions, we'll introduce a set of \textit{foundational functions} (\textit{F.F.}) which will be explored in \Cref{s:basic ff,s:advanced ff} and will be used in the rest of the paper.

\section{Basic foundational functions}\label{s:basic ff}
\subsection{Introduction}\label{ss:basic ff intro}

So far we have introduced the factorial and the Roman factorial. In this section we will set aside the factorial and its expansion to investigate a set of simple, Boolean-like functions that usually have binary outputs (0 or 1).

These so-called \textit{foundational functions} (\textit{F.F.}) are simple and are built upon each other, beginning with one called $\delta(n)$. They take integer or real values of $n$ and are categorized based on their outputs, which are similar to those of Boolean-valued functions. Basic information about Boolean algebra is provided in \hyperref[ss:add boolean algebra]{Addendum~\ref*{ss:add boolean algebra}}.

All function names have no correlation with other pre-existing functions of the same name, as all of them were named arbitrarily. For example, the function $\delta(n)$ used here has no correlation to the delta function involved in Fourier transformations.

Additionally, the variable $n$ will be used for the \textit{F.F.} instead of $x$ because the desired variable will usually take integer values and the letter $n$ is traditionally used for that purpose. On later parts of this study, $n$ will take real and complex values so the letter $z$ may be introduced.

Now, let's introduce the concept of output tables. Whenever a new function is defined in this paper, there will be a table showcasing its output values. For instance, the function $\sin\theta$ is represented by the following output table:
\begin{table}[H]
    \centering
    \begin{tabular}{c|*{5}{c}}
        $\theta$ & 0 & $\nicefrac{\pi}{2}$ & $\pi$ & $\nicefrac{3\pi}{2}$ & $2\pi$ \\
        \hline
        $\sin\theta$ & 0 & 1 & 0 & -1 & 0
    \end{tabular}
    \caption{The function $\sin\theta$}
    \label{tbl:sintheta}
\end{table}

Additionally, these output values will be expressed concisely as an "output pattern". \ref{tbl:sintheta} can be represented by the output pattern $[0, 1, 0, -1, 0]$. \textit{F.F.} in this paper will be relevant for a few values of $n$, so the output tables will have 3 or so columns.

The values of $n$ that will be used extensively are 0 and any positive or negative numbers. It should be noted that all \textit{F.F.} will converge at the limits as $n$ goes to infinity, unlike $\sin\theta$. Most of the time, \textit{foundational functions} will also have the same value across all positive or negative numbers.

Another example for an output table is this:
\begin{table}[H]
    \centering
    \begin{tabular}{c|*{5}{c}}
        $x$ & -$\infty$ & -1 & 0 & 1 & $+\infty$ \\
        \hline
        $e^{-x^2}$ & 0 & $\nicefrac{1}{e}$ & 1 & $\nicefrac{1}{e}$ & 0
    \end{tabular}
    \caption{The function $e^{-x^2}$}
    \label{tbl:e-x2}
\end{table}

This function can be expressed by the output pattern $[0, \nicefrac{1}{e}, 1, \nicefrac{1}{e}, 0]$. Most commonly though, output tables will have less columns. Using the previous example, we can condense \ref{tbl:e-x2} like so:
\begin{table}[H]
    \centering
    \begin{tabular}{c|*{3}{c}}
        $x$ & -$\infty$ & 0 & $+\infty$ \\
        \hline
        $e^{-x^2}$ & 0 & 1 & 0
    \end{tabular}
    \caption{The function $e^{-x^2}$, condensed}
    \label{tbl:e-x2 condensed}
\end{table}

In this output table, the values for $-\infty$ and $+\infty$ are determined by taking the limit of $e^{-x^2}$ as $x$ approaches $\spm$ infinity. However, for the functions presented in this section, there is no need to take limits because they will output a single value across their entire domain. The example shown above serves only to illustrate the concept of output tables.

Output tables are valuable tools for understanding how a function behaves across various input values. Throughout this paper, the input domain of foundational functions (\textit{F.F.}) will encompass all real numbers, even if they are evaluated at integer values.

In the rest of \Cref{s:basic ff}, as well as in \Cref{s:advanced ff}, we will define a set of functions denoted by (\textit{F.F.}). They will be helpful in \Cref{s:rf rec rl gen} where we will condense \ref{eq:rf pw} and \ref{eq:rf rec pw} into a single relationship each, and in \Cref{s:rf prod gen} in which we will build two universal product definitions for the Roman factorial. \textit{Foundational functions} will be necessary for intuitively comprehending the mathematical definitions presented in this paper.

\subsection{\texorpdfstring{The function $\delta(n)$}{The function delta of n}}\label{ss:delta function}

The first function to be defined here is named $\delta(n)$. For any integer $n$, $\delta(n)$ takes the form:
\begin{equation}\label{eq:simple delta}
    \delta(n) = n + \varepsilon \,, \M 0 < \varepsilon < 1 \,, \M n \in \mbZ \,.
\end{equation}

Choosing $\varepsilon = 0.5$, the resulting output table for this function is this:
\begin{center}
    \begin{tabular}{c|*{8}{c}}
        $n$ & $\ldots$ & -2 & -1 & 0 & 1 & 2 & $\ldots$ \\
        \hline
        $\delta(n)$ & $\ldots$ & -1.5 & -0.5 & 0.5 & 1.5 & 2.5 & $\ldots$
    \end{tabular}
\end{center}

The definition of $\delta(n)$ for integers is devised to ensure positive outputs for non-negative integers and negative outputs for negative integers. This requirement can be expressed through a piece-wise definition:
\begin{equation}\label{eq:delta pw}
    \delta(n) = 
    \begin{cases}
        + \,, \M n \geq 0 \\
        - \,, \M n < 0
    \end{cases}
    , \| n \in \mbZ \,.
\end{equation}

Here, the symbols $+$ or $-$ denote any positive or negative number, respectively. The function $\delta(n)$ aims to yield outputs with specific signs, the exact value of which is not relevant. Additionally, $\varepsilon$ cannot take the value of 0 or 1, because that would mean either $\delta(0) = 0$ or $\delta(-1) = 0$, and this is not desired.

Given that the precise output values are irrelevant, the previous output table can be simplified by indicating only the sign:
\begin{center}
    \begin{tabular}{c|*{9}{c}}
        $n$ & -$\infty$ & $\ldots$ & -2 & -1 & 0 & 1 & 2 & $\ldots$ & $\infty$ \\
        \hline
        $\delta(n)$ & $-$ & $\ldots$ & $-$ & $-$ & + & + & + & $\ldots$ & +
    \end{tabular}
\end{center}

Let's shorten the table to highlight the three key cases of interest: positive inputs, negative inputs and the zero input. The reduction results in this table:
\begin{table}[H]
    \centering
    \begin{tabular}{c|*{3}{c}}
        $n$ & $-$ & 0 & $+$ \\
        \hline
        $\delta(n)$ & $-$ & + & +
    \end{tabular}
    \caption{The function $\delta(n)$}
    \label{tbl:delta}
\end{table}

Hence, $\delta(n)$ exhibits the output pattern $[-, +, +]$. The parameter $\varepsilon$ can assume any value between 0 and 1 while adhering to the specified behavior.

This definition can be expanded to accept all real values of $n$ by utilizing the floor function:
\begin{equation}\label{eq:delta}
    \delta(n) = \ff{n} + \varepsilon \,, \M 0 < \varepsilon < 1 \,, \M n \in \mbR \,.
\end{equation}

The floor function, which effectively rounds down, is exemplified below:
\begin{gather*}
    \ff{1.3} = 1 \,, \qquad \ff{2.7} = 2 \,, \\[0.618034em]
    \ff{3} = 3 \,, \qquad \ff{-3.4} = -4 \,.
\end{gather*}

Additionally, the graph of the $\ff{n}$ is depicted next:
\begin{figure}[H]
    \centering
    \includegraphics[width = 1\linewidth]{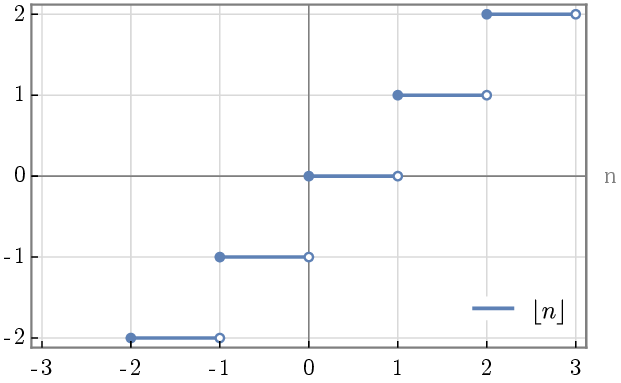}
    \caption{The floor function $\ff{n}$}
    \label{fig:floor}
\end{figure}

\ref{eq:delta} is termed the \textit{universal definition} of $\delta(n)$. Unlike \ref{eq:delta pw}, this formulation is not piece-wise and accommodates inputs from all real numbers.

Although the variable $n$ is initially limited to integer values, the floor function extends the definition of $\delta(n)$ across all real numbers. This definition also allows for plots devoid of discrete points, so that we can visualize clearly the behavior of a function.

Plotting this version of $\delta(n)$ with $\varepsilon = 0.5$ makes the following illustration:
\begin{figure}[H]
    \centering
    \includegraphics[width = 1\linewidth]{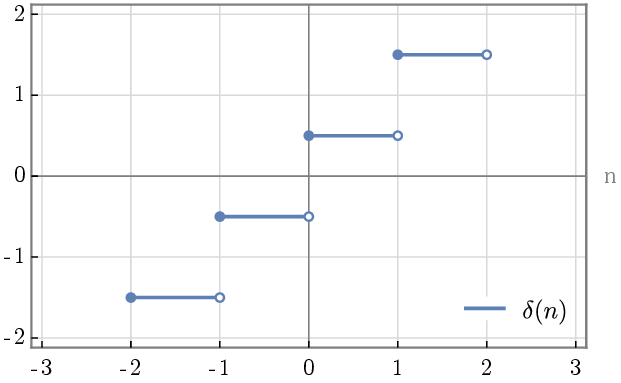}
    \caption{The function $\delta(n)$}
    \label{fig:delta floor}
\end{figure}

As observed, discontinuities arise naturally due to the floor function. Indeed, there is no continuous definition with $\delta(n)$'s outputs because of Bolzano's Theorem: for a continuous $\delta(n)$, it would be $\delta(n) = 0$ for some value between $\delta(-1) = -1 + \varepsilon < 0$ and $\delta(1) = \varepsilon > 0$, contradicting the requirement that $\delta(n)$ should not output zero for any input.

This function serves as a fundamental component for all subsequent \textit{F.F.} presented in this paper. A well-defined $\delta(n)$ for all real numbers is imperative\footnote{For a more advanced definition of $\delta(n)$ involving a Fourier approximation, see \hyperref[ss:add ff deltaf]{Addendum~\ref*{ss:add ff deltaf}}.}, even if non-integers are not utilized in this paper.

\subsection{\texorpdfstring{The function $\theta(n)$}{The function theta of n}}\label{ss:theta function}

The second \textit{F.F.}, denoted by $\theta(n)$, is simple. Its primary objective is to transform the positive and negative outputs of $\delta(n)$ into $+1$ or $-1$, respectively.

The piece-wise definition of $\theta(n)$ is as follows:
\begin{equation}\label{eq:theta pw}
    \theta(n) = 
    \begin{cases}
        +1 \,, \M n \geq 0 \\
        -1 \,, \M n < 0
    \end{cases}
    , \| n \in \mbR \,.
\end{equation}

This is achieved by dividing $\delta(n)$ by its absolute value. Since $\delta(n) \neq 0$, we can define $\theta(n)$ as such:
\begin{equation}\label{eq:theta}
    \theta(n) = \frac{\delta(n)}{\left|\delta(n)\right|} \:, \M n \in \mbR \,.
\end{equation}

\ref{eq:theta} represents the \textit{universal definition} of $\theta(n)$. Notably, the value of this function is consistently either 1 or -1, because of the division between identical terms. The sole variation is the sign of $\theta(n)$ which mirrors the sign of $\delta(n)$. The resulting output table is this:
\begin{table}[H]
    \centering
    \begin{tabular}{c|ccc}
        $n$ & $-$ & 0 & $+$ \\
        \hline
        $\theta(n)$ & -1 & 1 & 1
    \end{tabular}
    \caption{The function $\theta(n)$}
    \label{tbl:theta}
\end{table}

The graph of $\theta(n)$ is presented below:
\begin{figure}[H]
    \centering
    \includegraphics[width = 1\linewidth]{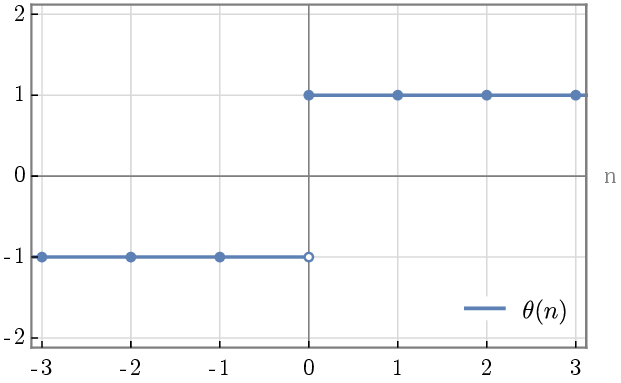}
    \caption{The function $\theta(n)$}
    \label{fig:theta}
\end{figure}

While resembling the $sign$ function's functionality, $\theta(n)$ is specifically chosen due to a critical distinction: $sign(0) = 0$, whereas $\theta(0) = 1$. This difference is important, and it is the rationale behind the definition of $\delta(n)$. A function similar to the $sign$ function, yet distinct from outputting zero, is necessary. The function $\theta(n)$ has the output pattern $[-1, 1, 1]$ and is crucial for defining \textit{F.F.} as well as building a universal Roman factorial in \Cref{s:rf rec rl gen}.

\subsection{\texorpdfstring{The function $\xi(n)$}{The function xi of n}}\label{ss:xi function}

The third \textit{F.F.} to be introduced here, denoted as $\xi(n)$, needs to have an output of 0 for all negative inputs and 1 for all others. In other words, the desirable output pattern is $[0, 1, 1]$. This can be concisely expressed through a closed-form definition:
\begin{equation}\label{eq:xi pw}
    \xi(n) = 
    \begin{cases}
        1 \,, \M n \geq 0 \\
        0 \,, \M n < 0
    \end{cases}
    , \| n \in \mbR \,.
\end{equation}

Let's construct $\xi(n)$ using $\theta(n)$ step by step. To begin with, we note that the output pattern of $\theta(n)$ is $[-1, 1, 1]$. Adding 1 to $\theta(n)$ essentially adds to each output, obtaining $[0, 2, 2]$. Dividing $\bp{1+ \theta(n)}$ by 2 yields the desired pattern $[0, 1, 1]$. Thus, we arrive at the \textit{universal definition} of $\xi(n)$:
\begin{equation}\label{eq:xi}
    \xi(n) = \frac{1 + \theta(n)}{2} \:, \M n \in \mbR \,.
\end{equation}

The plot of $\xi(n)$ is illustrated below:
\begin{figure}[H]
    \centering
    \includegraphics[width = 1\linewidth]{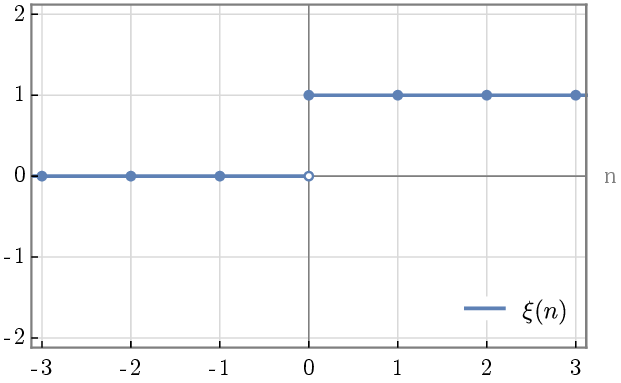}
    \caption{The function $\xi(n)$}
    \label{fig:xi}
\end{figure}

Correspondingly, the output table is presented as follows:
\begin{table}[H]
    \centering
    \begin{tabular}{c|ccc}
        $n$ & $-$ & 0 & $+$ \\
        \hline
        $\xi(n)$ & 0 & 1 & 1
    \end{tabular}
    \caption{The function $\xi(n)$}
    \label{tbl:xi}
\end{table}

This function has similarities to true/false states in Boolean logic, which is commonly used in programming. Specifically, the states $\xi(n \geq 0) = 1$ and $\xi(n < 0) = 0$ represent the true/false statements of the condition $n \geq 0$. Furthermore, it's noteworthy that the most subsequent \textit{F.F.} can also be expressed in Boolean algebra\footnote{More about Boolean algebra in \hyperref[ss:add boolean algebra]{Addendum~\ref*{ss:add boolean algebra}}.}.

\subsection{\texorpdfstring{The function $\xi^\prime(n)$}{The function xi prime of n}}\label{ss:xi prime function}

The function $\xi^\prime(n)$, marked with ($^\prime$), bears resemblance to $\xi(n)$, although it is not its derivative. In this paper, the ($^\prime$) symbol signifies similarity between functions rather than indicating a derivative. Notably, the \textit{F.F.} are devoid of derivatives since they are constructed upon the floor function, which itself lacks a derivative at integers. Derivatives will not be discussed within the scope of this work\footnote{It should be mentioned again that the names of the \textit{F.F.} are arbitrarily assigned Greek or Latin letters and have no relation to other functions with the same name. The naming process was conducted without any specific goal or consideration of potential overlaps with pre-existing functions.}.

The \textit{universal definition} of $\xi^\prime(n)$ is given as follows:
\begin{equation}\label{eq:xi prime}
    \xi^\prime(n) = \frac{1 - \theta(n)}{2} \:, \M n \in \mbR \,.
\end{equation}

Distinctively, the numerator has a subtraction instead of a summation, distinguishing it from $\xi(n)$. A relationship between $\xi(n)$ and $\xi^\prime(n)$ is observed:
\begin{equation}\label{eq:xi xi prime relation}
    \xi^\prime(n) = 1 - \xi(n) \,, \M n \in \mbR \,.
\end{equation}

In Boolean algebra, $\xi^\prime(n)$ yields the true/false states of the condition $n < 0 \,$, thereby exhibiting an output pattern of $[1, 0, 0]$.

The plot of $\xi^\prime(n)$ is depicted below:
\begin{figure}[H]
    \centering
    \includegraphics[width = 1\linewidth]{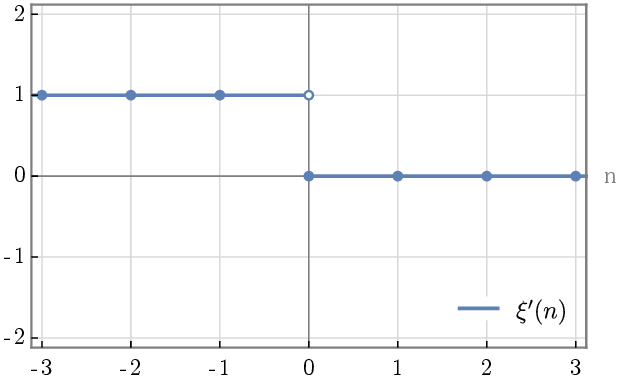}
    \caption{The function $\xi^\prime(n)$}
    \label{fig:xi prime}
\end{figure}

Thus far, four \textit{foundational functions} (\textit{F.F.}) have been defined: $\delta(n)$, $\theta(n)$, $\xi(n)$, and $\xi^\prime(n)$. Their respective outputs are summarized in \ref{tbl:basic ff summary no eta}:
\begin{table}[H]
    \centering
    \begin{tabular}{c|*{3}{c}}
        $n$ & $-$ & 0 & $+$ \\
        \hline
        $\delta(n)$ & $-$ & + & + \\
        $\theta(n)$ & -1 & 1 & 1 \\
        $\xi(n)$ & 0 & 1 & 1 \\
        $\xi^\prime(n)$ & 1 & 0 & 0
    \end{tabular}
    \caption{Outputs of $\delta(n)$, $\theta(n)$, $\xi(n)$ and $\xi^\prime(n)$}
    \label{tbl:basic ff summary no eta}
\end{table}

\subsection{\texorpdfstring{The function $\eta(n)$}{The function eta of n}}\label{ss:eta function}

The function $\eta(n)$ is defined by the following piece-wise definition:
\begin{equation}\label{eq:eta pw}
    \eta(n) = 
    \begin{cases}
        \,1 \,, & n \geq 0 \\
        \,1 \,, & n < 0 \,, \| n = \text{odd} \\
        -1 \,, & n < 0 \,, \| n = \text{even} \,.
    \end{cases}
\end{equation}

Unlike the previous functions, $\eta(n)$ varies based on whether $n$ is divisible by 2 or not. This function's ultimate purpose is to describe the alternating sign of the Roman factorial at negative integers.

Initially, we can address the sign-changing issue by considering the expression $(-1)^n$, which yields:
\begin{equation*}
    (-1)^n = 
    \begin{cases}
        -1 \,, & n = \text{odd} \\
        \,1 \,, & n = \text{even} \,.
    \end{cases}
\end{equation*}

This almost works, however the signs are opposite of what we specified in \ref{eq:eta pw}. This can be easily solved by raising $(-1)$ to the power $(n - 1)$:
\begin{equation*}
    (-1)^{n - 1} = 
    \begin{cases}
        \,1 \,, & n = \text{odd} \\
        -1 \,, & n = \text{even} \,.
    \end{cases}
\end{equation*}

Let's now address the case of positive inputs. Since we raise the base $(-1)$ to some power, if the base becomes $(+1)$ for $n > 0$, then the result will always be 1 regardless of the power. That is true because $1^n = 1, \| n \in \mbR \,$. To ensure that behavior we can use $\theta(n)$ for the base:
\begin{equation}\label{eq:eta wrong int}
    \eta(n) = \theta(n)^{n - 1} \,, \M n \in \mbZ \,.
\end{equation}

Also, by adjusting the exponent to $(\cf{n} - 1)$, where $\cf{n}$ is the ceiling function (rounding up), it is guaranteed that $(-1)$ is not raised to some fractional power for non-integer $n$. This initiative originates from the aim to define \textit{F.F.} for all real numbers, not just integers.

Now, we can express \ref{eq:eta pw} in this way:
\begin{equation}\label{eq:eta wrong}
    \eta(n) = \theta(n)^{\cf{n} - 1} \,, \M n \in \mbR \,.
\end{equation}

In the introduction to the Roman factorial, we observed a term similar to $\eta(n)$ in the piece-wise definition of \ref{eq:rf pw}. However, instead of having the power $(n - 1)$, it has $(-n - 1)$. Since we aim for $\eta(n)$ to be part of the universal Roman factorial definition, these terms should match exactly.

Fortunately, we can simply substitute $n$ with $-n$ directly within the power term of $\eta(n)$. This is because it doesn't matter if $(-1)$ is raised to a positive or negative power:
\begin{equation}\label{eq:minus 1 powers}
    (-1)^k \equiv (-1)^{-k} = \cfrac{1}{(-1)^k} = (-1)^k \,.
\end{equation}

For instance, the end result remains the same:
\begin{equation*}
    (-1)^{\cf{-5} - 1} = (-1)^{-\cf{-5} - 1} = 1 \,.
\end{equation*}

Ultimately, $\eta(n)$ is \textit{universally defined} in this way:
\begin{equation}\label{eq:eta}
    \eta(n) = \theta(n)^{-\cf{n} - 1} \,, \M n \in \mbR \,.
\end{equation}

The plot of $\eta(n)$ is shown below.
\begin{figure}[H]
    \centering
    \includegraphics[width = 1\linewidth]{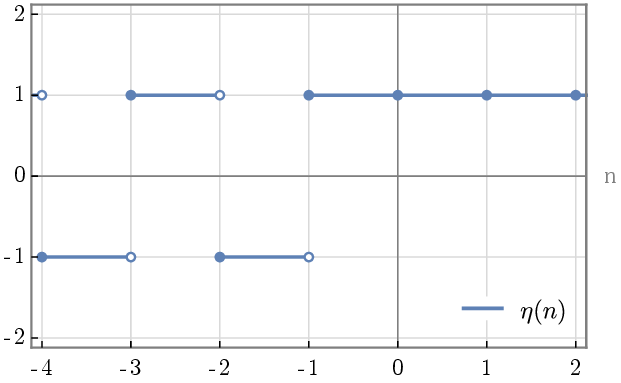}
    \caption{The function $\eta(n)$}
    \label{fig:eta}
\end{figure}

\subsection{Synopsis}\label{ss:synopsis of bff}

In this section, we have introduced various \textit{foundational functions} (\textit{F.F.}) that will play a crucial role throughout this paper. Given their similar naming conventions, it's easy to confuse them. Here's a summary of the five \textit{F.F.} covered thus far:
\begin{subequations}
    \begin{align}
        \delta(n) &= \ff{n} + 0.5 & [-, +, +] \tag{\ref{eq:delta}} \\[0.368034em]
        \theta(n) &= \cfrac{\delta(n)}{\left|\delta(n)\right|} & [-1, 1, 1] \tag{\ref{eq:theta}} \\[1.118034em]
        \xi(n) &= \cfrac{1 + \theta(n)}{2} & [0, 1, 1] \tag{\ref{eq:xi}} \\[0.868034em]
        \xi^\prime(n) &= \cfrac{1 - \theta(n)}{2} & [1, 0, 0] \tag{\ref{eq:xi prime}} \\[0.618034em]
        \eta(n) &= \theta(n)^{-\cf{n} - 1} & [\spm 1, 1, 1] \tag{\ref{eq:eta}}
    \end{align}
\end{subequations}

With these foundational functions established, we're now prepared to derive universal definitions for the Roman factorial.

You can refer to these equations for a reminder of the output patterns when these functions are utilized in \Cref{s:rf rec rl gen}. Additional \textit{F.F.} will be introduced in \Cref{s:advanced ff}, to be employed in \Cref{s:rf prod gen}.

\section{\texorpdfstring{Recursive generalizations of $\rf{n}$}{Recursive generalizations of the Roman factorial}}\label{s:rf rec rl gen}
\subsection{Introduction}\label{ss:rf rec rl gen intro}

In this section, we aim to unify the two parts of the piece-wise definition of the Roman factorial:
\begin{equation}\tag{\ref{eq:rf pw}}
    \rf{n} = 
    \begin{cases}
        n! & , \| n \in \mbZ^+_0 \\
        \cfrac{(-1)^{-n - 1}}{(-n - 1)!} & , \| n \in \mbZ^- ,
    \end{cases}
\end{equation}
in which the factorial is defined recursively as
\begin{equation}\tag{\ref{eq:tf rec}}
    n! = n(n - 1)! \,, \M 0! = 1 \,, \M n \in \mbZ^+ .
\end{equation}

Here, $\mbZ^+_0$ refers to the set of natural numbers plus 0, while $\mbZ^-$ represents negative integers. More about number sets in \hyperref[ss:add number sets]{Addendum~\ref*{ss:add number sets}}.

Additionally, we aim to rewrite the following closed-form definition as a singular expression as well:
\begin{equation}\tag{\ref{eq:rf rec pw}}
    \rf{n} = 
    \begin{cases}
        n \rf{n - 1} & , \| n \in \mbZ^+ \\[0.368034em]
        \cfrac{\rf{n + 1}}{n + 1} & , \| n \in \mbZ^- \setminus \{-1\} \,,
    \end{cases}
\end{equation}
where
\begin{equation}\label{eq:rf rec seeds}
    \rf{0} = \rf{-1} = 1 \,.
\end{equation}

We will call the generalization of \ref{eq:rf rec pw} as a \textit{recursive definition}, and that of \ref{eq:rf pw} as a \textit{Roman-like definition}. This is because we will address the easier one in the beginning and the more complicated one directly afterwards.

In order to achieve the unifications, we will need to modify step-by-step the two parts of each piece-wise definition so that they are exactly the same. Then the two cases will be identical and they will collapse into a function that works for all $n \in \mbZ\,$.

Despite the rewriting, the resulting relationships will still be based on a recursive factorial. Later, in \Cref{s:rf prod gen}, we will develop a new factorial definition from scratch.

\subsection{\texorpdfstring{Recursive definition step 1: $\theta(n)$}{Recursive definition step 1: theta of n}}\label{ss:rf rec gen step 1}

In \Cref{ss:rf rec gen step 1,ss:rf rec gen step 2} (the rest of this page) we aim to combine the two instances of the doubly-recursive definition of the Roman factorial:
\begin{equation}\tag{\ref{eq:rf rec pw}}
    \rf{n} = 
    \begin{cases}
        n \rf{n - 1} & , \| n \in \mbZ^+ \\[0.368034em]
        \cfrac{\rf{n + 1}}{n + 1} & , \| n \in \mbZ^- \setminus \{-1\} \,,
    \end{cases}
\end{equation}
where
\begin{equation}\tag{\ref{eq:rf rec seeds}}
    \rf{0} = \rf{-1} = 1 \,.
\end{equation}

Let's begin by revising the cases as follows:
\begin{equation*}
    \begin{cases}
        n \cdot \rf{n - 1} \\[0.368034em]
        \cfrac{1}{n + 1} \cdot \rf{n + 1}
    \end{cases}
    = 
    \begin{cases}
        n^{+1} \cdot \rf{n - (+1)} \\[0.368034em]
        (n + 1)^{-1} \cdot \rf{n - (-1)}
    \end{cases}
\end{equation*}

In this way, we can see that there are two occurrences of a function in each case. That function is $+1$ when $n \in \mbZ^+$ and $-1$ when $n \in \mbZ^- \setminus \{-1\}$ (meaning $n$ is a negative integer except $-1$).

We defined a function that meets these criteria. The function $\theta(n)$ has the form:
\begin{equation}\tag{\ref{eq:theta}}
    \theta(n) = \cfrac{\delta(n)}{|\delta(n)|} =
    \begin{cases}
        1 & , \| n \in \mbZ^+ \\ 
        -1 & , \| n \in \mbZ^- \setminus \{-1\} \,.
    \end{cases}
\end{equation}

We can replace those occurrences with the \textit{F.F.} $\theta(n)$, since it outputs the same values at the domains of interest:
\begin{equation}\label{eq:rf rec gen step 1.0}
    \begin{cases}
        n^{\theta(n)} \cdot \rf{n - \theta(n)} \\[0.368034em]
        (n + 1)^{\theta(n)} \cdot \rf{n - \theta(n)}
    \end{cases}
\end{equation}

Now, our definition has the following form:
\begin{equation}\label{eq:rf rec gen step 1.1}
    \rf{n} = 
    \begin{cases}
        n^{\theta(n)} \rf{n - \theta(n)} & , \| n \in \mbZ^+ \\[0.368034em]
        (n + 1)^{\theta(n)} \rf{n - \theta(n)} & , \| n \in \mbZ^- \setminus \{-1\} \,,
    \end{cases}
\end{equation}
where
\begin{equation}\tag{\ref{eq:rf rec seeds}}
    \rf{0} = \rf{-1} = 1 \,.
\end{equation}

\subsection{\texorpdfstring{Recursive definition step 2: $\xi^\prime(n)$}{Recursive definition step 2: xi prime of n}}\label{ss:rf rec gen step 2}

We observe that in the \ref{eq:rf rec gen step 1.0}, the only difference between the cases is the addition of 1 in the first term, which is absent when $n \in \mbZ^+$.

Let's highlight this absence like so:
\begin{equation}\label{eq:rf rec gen step 2.0}
    \rf{n} = 
    \begin{cases}
        (n + 0)^{\theta(n)} \rf{n - \theta(n)} & , \| n \in \mbZ^+ \\[0.368034em]
        (n + 1)^{\theta(n)} \rf{n - \theta(n)} & , \| n \in \mbZ^- \setminus \{-1\} \,.
    \end{cases}
\end{equation}

Notice that this additive term is 0 when $n \in \mbZ^+$ and 1 when $n \in \mbZ^- \setminus \{-1\}$. This behavior is described by the \textit{F.F.} $\xi^\prime(n)$:
\begin{equation}\tag{\ref{eq:xi prime}}
    \xi^\prime(n) = \cfrac{1 - \theta(n)}{2} =
    \begin{cases}
        0 & , \| n \in \mbZ^+ \\
        1 & , \| n \in \mbZ^- \setminus \{-1\} \,.
    \end{cases}
\end{equation}

Incorporating this \textit{F.F.} into \ref{eq:rf rec gen step 2.0} results in the following formulation:
\begin{equation*}
    \begin{cases}
        \bp{n + \xi^\prime(n)}^{\theta(n)} \rf{n - \theta(n)} & , \| n \in \mbZ^+ \\[0.368034em]
        \bp{n + \xi^\prime(n)}^{\theta(n)} \rf{n - \theta(n)} & , \| n \in \mbZ^- \setminus \{-1\}
    \end{cases}
\end{equation*}

Thus, since the two cases have exactly the same form, we can consolidate them into one:
\begin{equation}\label{eq:rf rec gen}
    \rf{n} = \bp{n + \xi^\prime(n)}^{\theta(n)} \rf{n - \theta(n)} \,, \M n \in \mbZ \setminus \{0, \|-1\} \,,
\end{equation}
where
\begin{equation}\tag{\ref{eq:rf rec seeds}}
    \rf{0} = \rf{-1} = 1 \,.
\end{equation}

This result concludes the recursive generalization, so in the next subsections we will address the Roman-like one.

\subsection{\texorpdfstring{Roman-like definition step 1: $\eta(n)$}{Roman-like definition step 1: eta of n}}\label{ss:rf gen step 1}

In \Cref{ss:rf gen step 1,ss:rf gen step 2,ss:rf gen step 3,ss:rf gen step 4} we aim to unite the two cases of the original definition of the Roman factorial, expressed as
\begin{equation}\tag{\ref{eq:rf pw}}
    \rf{n} = 
    \begin{cases}
        n! & , \| n \in \mbZ^+_0 \\
        \cfrac{(-1)^{-n - 1}}{(-n - 1)!} & , \| n \in \mbZ^- ,
    \end{cases}
\end{equation}
in which the factorial is defined recursively as
\begin{equation}\tag{\ref{eq:tf rec}}
    n! = n(n - 1)! \,, \M 0! = 1 \,, \M n \in \mbZ^+ .
\end{equation}

To start, let's separate the numerator from the fraction in the second case. Additionally, it's important to emphasize that the factorial in the first case is multiplied by 1:
\begin{equation}\label{eq:rf gen step 1.0}
    \rf{n} = 
    \begin{cases}
        1 \cdot n! & , \| n \in \mbZ^+_0 \\
        (-1)^{-n - 1} \cdot \cfrac{1}{(-n - 1)!} & , \| n \in \mbZ^- .
    \end{cases}
\end{equation}

In this manner, it's evident that the isolated term corresponds to the \textit{F.F.} $\eta(n)$. Below, we showcase the function as it is defined for integers:
\begin{equation}\tag{\ref{eq:eta}}
    \eta(n) = \theta(n)^{-n - 1} =
    \begin{cases}
        1 & , \| n \in \mbZ^+_0 \\
        (-1)^{-n - 1} & , \| n \in \mbZ^- .
    \end{cases}
\end{equation}

Even though we initially defined $\eta(n)$ using the floor function for $n \in \mbR\,$, it's important to note that the Roman factorial is defined solely for integer values of $n$, which is a subset of $\mbR\,$. Thus, we can safely substitute it within the definition:
\begin{equation}\label{eq:rf gen step 1.1}
    \rf{n} = 
    \begin{cases}
        \eta(n) \cdot n! & , \| n \in \mbZ^+_0 \\
        \eta(n) \cdot \cfrac{1}{(-n - 1)!} & , \| n \in \mbZ^- .
    \end{cases}
\end{equation}

Understanding this procedure of substituting \textit{F.F.} step-by-step is crucial in understanding the \textit{generalization process} presented in this paper. The next subsections will follow a similar approach with different \textit{F.F.} focused on each step. This method will be pivotal for constructing definitions all throughout this paper.

\subsection{\texorpdfstring{Roman-like definition step 2: $\theta(n)$}{Roman-like definition step 2: theta of n}}\label{ss:rf gen step 2}

For the second step, let's observe that the factorial in the second case is inverted, while the factorial in the first case remains unchanged. In other words, the factorial of $n$ or $(-n - 1)$ is raised to the power of 1 or $-1$, respectively:
\begin{equation}\label{eq:rf gen step 2.0}
    \rf{n} = 
    \begin{cases}
        \eta(n) \cdot \bp{n!}^1 & , \| n \in \mbZ^+_0 \\[0.368034em]
        \eta(n) \cdot \bb{(-n - 1)!}^{-1} & , \| n \in \mbZ^- .
    \end{cases}
\end{equation}

To express this change, we need a function for the exponent that equals 1 when $n$ is non-negative and -1 when $n$ is negative. Another way to describe this requirement is the output pattern $[-1, 1, 1]$.

The \textit{F.F.} that satisfies this demand is $\theta(n)$. Its definition is shown below\footnote{We note that all \textit{foundational functions} (\textit{F.F.}) were defined for a broader set than what is required in this section. The definitions of \textit{F.F.} rewritten here are specifically for integers, although they are well defined for all real numbers in \Cref{s:basic ff}.}:
\begin{equation}\tag{\ref{eq:theta}}
    \theta(n) = \cfrac{\delta(n)}{|\delta(n)|} =
    \begin{cases}
        1 & , \| n \in \mbZ^+_0 \\ 
        -1 & , \| n \in \mbZ^- .
    \end{cases}
\end{equation}

Let's incorporate $\theta(n)$ into our revised Roman factorial definition:
\begin{equation}\label{eq:rf gen step 2.1}
    \rf{n} = 
    \begin{cases}
        \eta(n) \cdot \bp{n!}^{\,\theta(n)} & , \| n \in \mbZ^+_0 \\[0.368034em]
        \eta(n) \cdot \bb{(-n - 1)!}^{\,\theta(n)} & , \| n \in \mbZ^- .
    \end{cases}
\end{equation}

We will omit the parentheses enclosing the factorial. Now, it is more concise and clear to the eye:
\begin{equation}\label{eq:rf gen step 2.2}
    \rf{n} = 
    \begin{cases}
        \eta(n) \cdot n!^{\,\theta(n)} & , \| n \in \mbZ^+_0 \\[0.368034em]
        \eta(n) \cdot (-n - 1)!^{\,\theta(n)} & , \| n \in \mbZ^- .
    \end{cases}
\end{equation}

So far, we have introduced two \textit{foundational functions} into \ref{eq:rf pw} in an attempt to rewrite it compactly. Before introducing another \textit{F.F.} however, there is one more adjustment we can make involving the absolute value of $n$.

\subsection{\texorpdfstring{Roman-like definition step 3: $|n|$}{Roman-like definition step 3: Absolute value of n}}\label{ss:rf gen step 3}

In the second case of \ref{eq:rf gen step 2.2}, there is a negative sign before the $n$ in the term $(-n - 1)$. This ensures that the factorial calculated is always positive, because if $n < 0$, then $-n > 0$. In contrast, in the first case, the $n$ inside the factorial is written without any sign (implying the positive sign).

This behavior can be described by the absolute value of $n$, which is defined like this:
\begin{equation}\label{eq:abs}
    |n| = 
    \begin{cases}
        n & , \| n \geq 0 \\
        -n & , \| n < 0 \,.
    \end{cases}
\end{equation}

We can substitute $|n|$ into \ref{eq:rf gen step 2.2} without any further consideration, since the absolute value is a function defined for all reals and thus all integers.

Now, the piece-wise definition of the Roman factorial can be expressed like this:
\begin{equation}\label{eq:rf gen step 3.0}
    \rf{n} = 
    \begin{cases}
        \eta(n) \cdot |n|!^{\,\theta(n)} & , \| n \in \mbZ^+_0 \\[0.368034em]
        \eta(n) \cdot \bp{|n| - 1}!^{\,\theta(n)} & , \| n \in \mbZ^- .
    \end{cases}
\end{equation}

We are nearly at the point where both cases have identical expressions. The only remaining difference is the $-1$ inside the factorial of the second case, which we'll address next.

\subsection{\texorpdfstring{Roman-like definition step 4: $\xi^\prime(n)$}{Roman-like definition step 4: xi prime of n}}\label{ss:rf gen step 4}

Let's investigate the $-1$ that appears only when $n \in \mbZ^-$. We can begin by expanding the first case as shown below:
\begin{equation}\label{eq:rf gen step 4.0}
    \rf{n} = 
    \begin{cases}
        \eta(n) \cdot \bp{|n| - 0}!^{\,\theta(n)} & , \| n \in \mbZ^+_0 \\[0.368034em]
        \eta(n) \cdot \bp{|n| - 1}!^{\,\theta(n)} & , \| n \in \mbZ^- .
    \end{cases}
\end{equation}

It is evident that we need a function to subtract from $|n|$, which equals 0 for $n \in \mbZ^+_0$ and 1 for $n \in \mbZ^-$. Alternatively, we desire a function with the output pattern $[1, 0, 0]$. This pattern adheres to the \textit{F.F.} $\xi^\prime(n)$:
\begin{equation}\tag{\ref{eq:xi prime}}
    \xi^\prime(n) = \frac{1 - \theta(n)}{2} = 
    \begin{cases}
        0 & , \| n \in \mbZ^+_0 \\
        1 & , \| n \in \mbZ^- .
    \end{cases}
\end{equation}

Thus, the two parts can become exactly the same:
\begin{equation}\label{eq:rf gen step 4.1}
    \rf{n} = 
    \begin{cases}
        \eta(n) \cdot \bp{|n| - \xi^\prime(n)}!^{\,\theta(n)} & , \| n \in \mbZ^+_0 \\[0.368034em]
        \eta(n) \cdot \bp{|n| - \xi^\prime(n)}!^{\,\theta(n)} & , \| n \in \mbZ^- .
    \end{cases}
\end{equation}

Finally, we have reached the recursive universal definition of integer Roman factorials. It is written in terms of the \textit{F.F.} as follows:
\begin{equation}\label{eq:rf gen}
    \rf{n} = \eta(n) \cdot \bp{|n| - \xi^\prime(n)}!^{\,\theta(n)} \,, \M n \in \mbZ \,,
\end{equation}
where
\begin{equation}\tag{\ref{eq:tf rec}}
    n! = n(n - 1)! \,, \M 0! = 1 \,, \M n \in \mbZ^+ .
\end{equation}

This marks the completion of one of our objectives. To sum up, we established a collection of \textit{foundational functions}, which enabled us to express the definition of the Roman factorial concisely and efficiently. In the following section, we will introduce another set of \textit{F.F.} which will be used alongside those defined in \Cref{s:basic ff} in order to build advanced factorial definitions. For example, in \Cref{s:rf prod gen}, we aim to create a factorial similar to the Roman factorial, but expressed as a non-recursive product.

\section{Advanced foundational functions}\label{s:advanced ff}
\subsection{Introduction}\label{ss:advanced ff intro}

In this section we will present 5 more \textit{foundational functions}. These will be crafted step-by-step using the other functions presented so far, and they will all be necessary in our endeavor to construct a non-recursive definition for the Roman factorial.

\subsection{\texorpdfstring{The function $\Theta(n)$}{The function Theta of n}}\label{ss:Theta function}

The first advanced function\footnote{These \textit{advanced} functions are not more difficult to understand, but they are more intricate, as they consist of many previously defined \textit{foundational functions}.} introduced in this section, denoted by $\Theta(n)$, is desired to output 1 only when $n = 0$ and 0 in all other cases. This can be expressed in a piece-wise definition like so:
\begin{equation}\label{eq:Theta pw}
    \Theta(n) = 
    \begin{cases}
        1 \,, \M n = 0 \\
        0 \,, \M n \neq 0
    \end{cases}
    , \| n \in \mbR \,.
\end{equation}

We aim to establish a universal definition for this function, similar to the ones we've devised previously.

Let's start by examining the function $\xi(n)$. We know it follows the output pattern $[0, 1, 1]$, which is close to the desired output pattern for $\Theta(n)$: $[0, 1, 0]$.

One approach to express $\Theta(n)$ is to multiply $\xi(n)$ by a function described by the outputs $[1, 1, 0]$. The function $\xi^\prime(n)$ is close, with the pattern $[1, 0, 0]$. However, $\xi^\prime(0) = 0$, whereas we want a function that equals 1 for $n = 0$.

Fortunately, there is a trick we can use to find the desired function. When we change the sign of the input $n$ in $\xi(n)$, the output pattern also reverses to $[1, 1, 0]$. Thus, the function we seek is $\xi(-n)$.

Therefore, we define $\Theta(n)$ \textit{universally} as:
\begin{equation}\label{eq:Theta}
    \Theta(n) = \xi(n) \cdot \xi(-n) \,, \M n \in \mbR \,.
\end{equation}

The next figure displays the graph of $\Theta(n)$:
\begin{figure}[H]
    \centering
    \includegraphics[width = 1\linewidth]{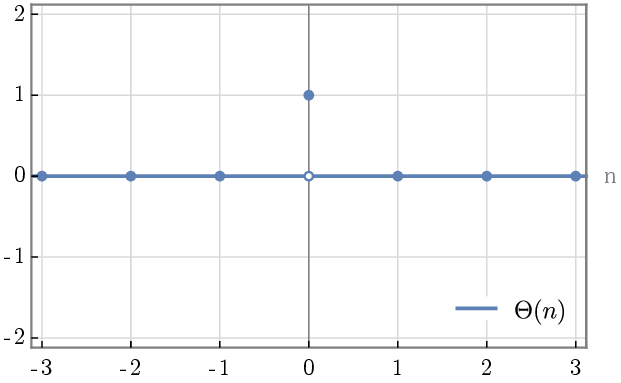}
    \caption{The function $\Theta(n)$}
    \label{fig:Theta}
\end{figure}

We also present the output table of $\Theta(n)$, as we did previously with \textit{F.F.}:
\begin{table}[H]
    \centering
    \begin{tabular}{c|*{3}{c}}
        $n$ & $-$ & 0 & $+$ \\
        \hline
        $\Theta(n)$ & 0 & 1 & 0
    \end{tabular}
    \caption{The function $\Theta(n)$}
    \label{tbl:Theta}
\end{table}

In contrast to previous \textit{foundational functions}, $\Theta(n)$ appears as a straight line on its plot, except for one point that stands out. In terms of Boolean algebra, $\Theta(n)$ represents the condition of whether $n$ equals zero.

$\Theta(n)$ can also be derived from the basic \textit{F.F.} in more than one way. Here is a list of some of them:
\begin{equation}\label{eq:Theta alts}
    \Theta(n) = 
    \begin{cases}
        \cfrac{\theta(n) + \theta(-n)}{2} \\[1.118034em]
        \cfrac{\xi(n) + \xi(-n)}{2} \\[1.118034em]
        \cfrac{1 + \theta(n) \cdot \theta(-n)}{2} \\[1.118034em]
        \cfrac{1 + \theta\bp{-|n|}}{2} \\[0.868034em]
        \xi\bp{-|n|}
    \end{cases}
    , \| n \in \mbR \,.
\end{equation}

Further exploration into alternative definitions of $\Theta(n)$ is possible, although it's not necessary to exhaustively enumerate every possibility. The alternatives provided in \ref{eq:Theta alts} represent a subset discovered during experimentation with the introduced functions thus far.

As an exercise, one can delve deeper into understanding these relationships by proving them. This can be accomplished through mathematical proofs or by visualizing the individual functions' graphs\footnote{All graphs in this paper were made with Wolfram Mathematica. For help with the relevant code, contact information is at \Cref{s:references}.} and their combinations.

\subsection{\texorpdfstring{The function $Q(n)$}{The function Q of n}}\label{ss:Q function}

The function $Q(n)$ represents the $sign$ function, which exhibits the output pattern $[-1, 0, 1]$. It returns the sign of a number, like $\theta(n)$, but it equals 0 if that number is 0. Its piece-wise definition is this:
\begin{equation}\label{eq:Q pw}
    Q(n) = 
    \begin{cases}
        \,1 \,, & n > 0 \\
        \,0 \,, & n = 0 \\
        -1 \,, & n < 0
    \end{cases}
    \M , \| n \in \mbR \,.
\end{equation}

Let's attempt to unify the parts of this definition. We can start by observing that $\theta(n)$ is almost identical to $Q(n)$, with the difference that $Q(0) = 0$ whereas $\theta(0) = 1$. This suggests that we can introduce $\theta(n)$ in the cases where $Q(n) = \theta(n)$:
\begin{equation}\label{eq:Q pw step 0}
    Q(n) = 
    \begin{cases}
        \theta(n) \,, & n > 0 \\
        0 \,, & n = 0 \\
        \theta(n) \,, & n < 0
    \end{cases}
    \M , \| n \in \mbR \,.
\end{equation}

Since the first and third case have identical expressions, we can combine them into one by merging their domains to $n \neq 0$. Thus, we have reduced the number of cases from three down to two:
\begin{equation}\label{eq:Q pw step 1}
    Q(n) = 
    \begin{cases}
        \theta(n) \,, & n \neq 0 \\[0.118034em]
        0 \,, & n = 0
    \end{cases}
    \M , \| n \in \mbR \,.
\end{equation}

Next, we will attempt to introduce another \textit{F.F.} in the definition of $Q(n)$. Let's start by examining the second case: we know that $0 = \theta(0) - 1$, so we can substitute that expression in the definition. Additionally, let's rewrite the first case as $\theta(n) - 0$, in order to highlight a subtraction from $\theta(n)$:
\begin{equation}\label{eq:Q pw step 2}
    Q(n) = 
    \begin{cases}
        \theta(n) - 0 \,, & n \neq 0 \\[0.118034em]
        \theta(n) - 1 \,, & n = 0
    \end{cases}
    \M , \| n \in \mbR \,.
\end{equation}

The two cases look very similar, the only thing that remains is to find a function to subtract from $\theta(n)$. Since we know that $\Theta(n)$ has the output pattern $[0, 1, 0]$, we can subtract it from $\theta(n)$. Thus, the \textit{universal definition} of $Q(n)$ is as follows:
\begin{equation}\label{eq:Q}
    Q(n) = \theta(n) - \Theta(n) \,, \M n \in \mbR \,.
\end{equation}

The plot of $Q(n)$ is depicted below, along with the corresponding output table:
\begin{figure}[H]
    \centering
    \includegraphics[width = 1\linewidth]{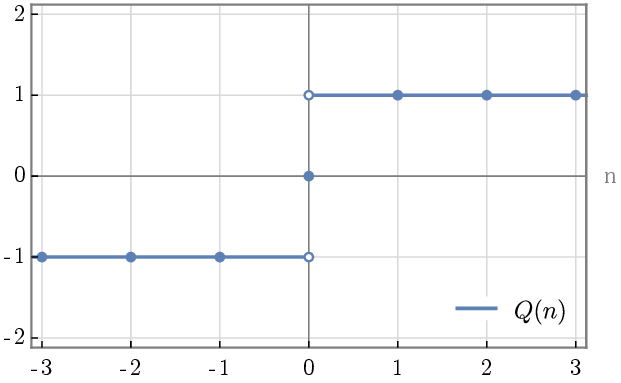}
    \caption{The function $Q(n)$}
    \label{fig:Q}
\end{figure}

\begin{table}[H]
    \centering
    \begin{tabular}{c|*{3}{c}}
        $n$ & $-$ & 0 & $+$ \\
        \hline
        $Q(n)$ & -1 & 0 & 1
    \end{tabular}
    \caption{The function $Q(n)$}
    \label{tbl:Q}
\end{table}

This function is useful because it behaves identically to $sign(n)$, producing $+1$ or -$1$ for positive or negative inputs respectively, and 0 for the input 0. However, $Q(n)$ is defined using \textit{foundational functions}, and henceforth, it will be used in place of $sign(n)$.

Some alternative versions of $Q(n)$ are listed here:
\begin{equation}\label{eq:Q alts}
    Q(n) = 
    \begin{cases}
        \cfrac{\theta(n) - \theta(-n)}{2} \\[0.868034em]
        \xi(n) - \xi(-n) \\[0.618034em]
        \xi^\prime(-n) - \xi^\prime(n)
    \end{cases}
    , \| n \in \mbR \,.
\end{equation}

\subsection{\texorpdfstring{The function $Q^\prime(n)$}{The function Q prime of n}}\label{ss:Q prime function}

The function $Q^\prime(n)$ is akin to $\Theta(n)$, but with reversed criteria. In essence, $Q^\prime(n)$ yields 0 only when $n = 0$ and 1 for $n \neq 0$. Thus, the piece-wise definition of $Q^\prime(n)$ is as follows:
\begin{equation}\label{eq:Q prime pw}
    Q^\prime(n) = 
    \begin{cases}
        0 \,, & n = 0 \\
        1 \,, & n \neq 0
    \end{cases}
    \M , \| n \in \mbR \,.
\end{equation}

An easy approach to achieve this goal using the \textit{F.F.} is to take the absolute value of $Q(n)$. Since $Q(n)$ has the output pattern $[-1, 0, 1]$, we see that $|Q(n)|$ will have the pattern $[1, 0, 1]$. That is our desired set of outputs, so we can define $Q^\prime(n)$ as the absolute value of $Q(n)$:
\begin{equation}\label{eq:Q prime}
    Q^\prime(n) = |Q(n)| \,, \M n \in \mbR \,.
\end{equation}

The output table for this function is provided below, as well as its plot:
\begin{table}[H]
    \centering
    \begin{tabular}{c|*{3}{c}}
        $n$ & $-$ & 0 & $+$ \\
        \hline
        $Q^\prime(n)$ & 1 & 0 & 1
    \end{tabular}
    \caption{The function $Q^\prime(n)$}
    \label{tbl:Q prime}
\end{table}

\begin{figure}[H]
    \centering
    \includegraphics[width = 1\linewidth]{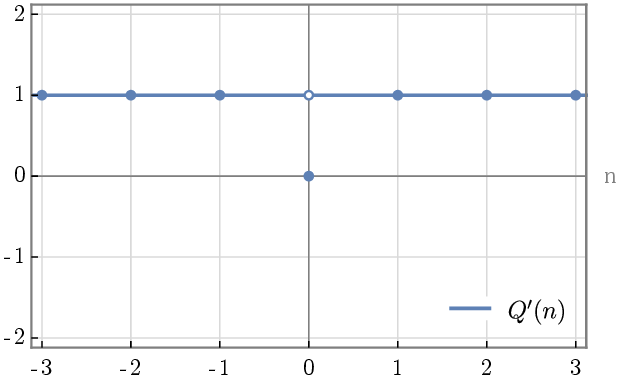}
    \caption{The function $Q^\prime(n)$}
    \label{fig:Q prime}
\end{figure}

Another approach is to subtract $\Theta(n)$ from 1, which is equally valid. Generally, the Boolean-like \textit{foundational functions} presented in \Cref{s:advanced ff} can be combined in many ways to achieve a certain output pattern. Below you will find more alternative definitions of this function:
\begin{equation}\label{eq:Q prime alts}
    Q^\prime(n) = 
    \begin{cases}
        1 - \Theta(n) \\[0.868034em]
        \cfrac{1 - \theta(n) \cdot \theta(-n)}{2} \\[0.868034em]
        \xi^\prime\bp{-|n|}
    \end{cases}
    , \| n \in \mbR \,.
\end{equation}

\subsection{\texorpdfstring{The function $\Psi(n)$}{The function Psi of n}}\label{ss:Psi function}

The next function we'll introduce in this paper is $\Psi(n)$. We need this function to be equal to $n$ for all $n \neq 0$, and 1 for $n = 0$. Its requirements can be expressed in the following piece-wise definition:
\begin{equation}\label{eq:Psi pw}
    \Psi(n) = 
    \begin{cases}
        n \,, & n \neq 0 \\
        1 \,, & n = 0
    \end{cases}
    \M , \| n \in \mbR \,.
\end{equation}

Defining this function is actually very straightforward. We can simply add $\Theta(n)$ to $n$, like so:
\begin{equation}\label{eq:Psi}
    \Psi(n) = n + \Theta(n) \,, \M n \in \mbR \,.
\end{equation}

This \textit{F.F.} is useful in many cases, for instance when we want to avoid division by 0. Dividing by $\Psi(n)$ ensures that there will never be any anomalies.

The outputs of $\Psi(n)$ and its plot are listed below:
\begin{table}[H]
    \centering
    \begin{tabular}{c|*{3}{c}}
        $n$ & $-$ & 0 & $+$ \\
        \hline
        $\Psi(n)$ & $n$ & 1 & $n$
    \end{tabular}
    \caption{The function $\Psi(n)$}
    \label{tbl:Psi}
\end{table}

\begin{figure}[H]
    \centering
    \includegraphics[width = 1\linewidth]{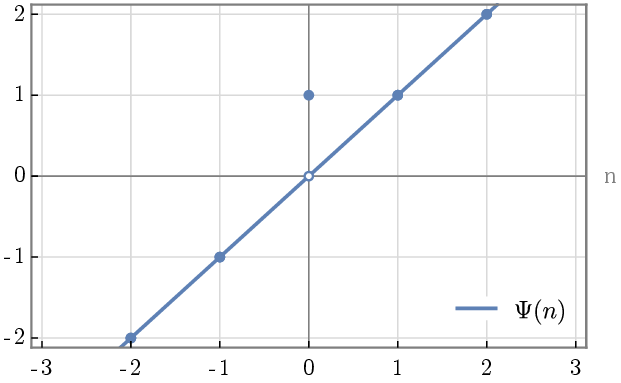}
    \caption{The function $\Psi(n)$}
    \label{fig:Psi}
\end{figure}

Additionally, we can interpret $\Psi(n)$ as a function that maps $0$ to $1$ while leaving all other values of $n$ unchanged. The simplicity of this function makes it useful for various mathematical operations, but it's most useful in defining the next \textit{F.F.}, $\Phi(n)$.

\subsection{\texorpdfstring{The function $\Phi(n)$}{The function Phi of n}}\label{ss:Phi function}

The final function to be presented in this paper is $\Phi(n)$. This function has requirements that can be expressed in the following piece-wise definition:
\begin{equation}\label{eq:Phi pw}
    \Phi(n) = 
    \begin{cases}
        1 \,, & n \geq 0 \\
        n \,, & n < 0
    \end{cases}
    \M , \| n \in \mbR \,.
\end{equation}

To meet these requirements, we can use an idea similar to the one employed in the definition of $\eta(n)$. In that case, we wanted to satisfy the desired outputs only for negative $n$, so we raised $\theta(n)$ to a power. Similarly, here we want $\Phi(n)$ to be 1 when $n > 0$. We can achieve this by raising $n$ to a power that is 0 when $n \geq 0$ and 1 when $n < 0$. This power is the function $\xi^\prime(n)$:
\begin{equation}\label{eq:Phi wrong}
    \Phi(n) = n^{\xi^\prime(n)} \,, \M n \in \mbR \,.
\end{equation}

This works well, except for the case when $n = 0$: $\Phi(0) = 0^0$, which is problematic because $0^0$ is undefined\footnote{The indeterminate form $0^0$ does not have a universally agreed-upon value, although it is given the value 1 in a lot of contexts \cite{0tothe0}. More about $0^0$ in \hyperref[ss:add 0 to the 0]{Addendum~\ref*{ss:add 0 to the 0}}.}. Fortunately we can solve this problem by using $\Psi(n)$. We end up with this \textit{universal definition} of $\Phi(n)$:
\begin{equation}\label{eq:Phi}
    \Phi(n) = {\Psi(n)}^{\xi^\prime(n)} \,, \M n \in \mbR \,.
\end{equation}

The output table and plot of $\Phi(n)$ are shown below:
\begin{table}[H]
    \centering
    \begin{tabular}{c|*{3}{c}}
        $n$ & $-$ & 0 & $+$ \\
        \hline
        $\Phi(n)$ & $n$ & 1 & 1
    \end{tabular}
    \caption{The function $\Phi(n)$}
    \label{tbl:Phi}
\end{table}

\begin{figure}[H]
    \centering
    \includegraphics[width = 1\linewidth]{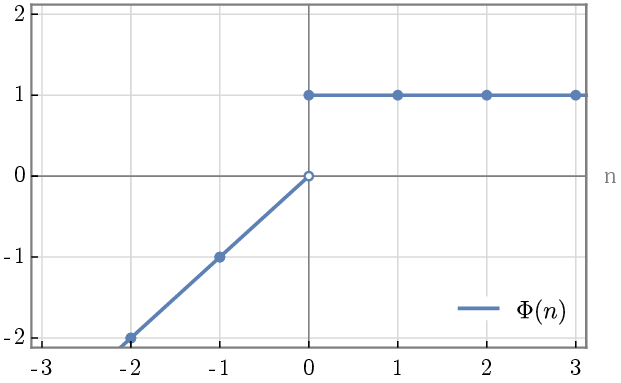}
    \caption{The function $\Phi(n)$}
    \label{fig:Phi}
\end{figure}

\subsection{Synopsis}\label{ss:synopsis of aff}

In short, this section is about presenting another set of 5 \textit{foundational functions}. They were all expressed in terms of other \textit{F.F.}, and along with those in \Cref{s:basic ff} they constitute a set of Boolean-like functions that are instrumental to this paper.

In the following table, we list every \textit{foundational function} defined in \Cref{s:basic ff,s:advanced ff} for $n \in \mbR$:
\begin{align}
    \delta(n) & = \ff{n} + 0.5 & [-, +, +] \tag{\ref{eq:delta}} \\[0.368034em]
    \theta(n) & = \cfrac{\delta(n)}{\left|\delta(n)\right|} & [-1, 1, 1] \tag{\ref{eq:theta}} \\[0.618034em]
    \xi(n) & = \cfrac{1 + \theta(n)}{2} & [0, 1, 1] \tag{\ref{eq:xi}} \\[0.618034em]
    \xi^\prime(n) & = \cfrac{1 - \theta(n)}{2} & [1, 0, 0] \tag{\ref{eq:xi prime}} \\[0.618034em]
    \eta(n) & = \theta(n)^{-\cf{n} - 1} & [\spm 1, 1, 1] \tag{\ref{eq:eta}} \\[0.368034em]
    \Theta(n) & = \xi(n) \cdot \xi(-n) & [0, 1, 0] \tag{\ref{eq:Theta}} \\[0.368034em]
    Q(n) & = \theta(n) - \Theta(n) & [-1, 0, 1] \tag{\ref{eq:Q}} \\[0.368034em]
    Q^\prime(n) & = 1 - \Theta(n) & [1, 0, 1] \tag{\ref{eq:Q prime}} \\[0.368034em]
    \Psi(n) & = n + \Theta(n) & [n, 1, n] \tag{\ref{eq:Psi}} \\[0.368034em]
    \Phi(n) & = {\Psi(n)}^{\xi^\prime(n)} & [n, 1, 1] \tag{\ref{eq:Phi}}
\end{align}

Also, we include an output table of all \textit{F.F.}:
\begin{table}[H]
    \centering
    \begin{tabular}{c|*{3}{c}}
        $n$ & $-$ & 0 & $+$ \\
        \hline
        $\delta(n)$ & $-$ & + & + \\
        $\theta(n)$ & -1 & 1 & 1 \\
        $\xi(n)$ & 0 & 1 & 1 \\
        $\xi^\prime(n)$ & 1 & 0 & 0 \\
        $\eta(n)$ & $\spm 1$ & 1 & 1 \\
        $\Theta(n)$ & 0 & 1 & 0 \\
        $Q(n)$ & -1 & 0 & 1 \\
        $Q^\prime(n)$ & 1 & 0 & 1 \\
        $\Psi(n)$ & $n$ & 1 & $n$ \\
        $\Phi(n)$ & $n$ & 1 & 1
    \end{tabular}
    \caption{Synopsis of the \textit{F.F.}}
    \label{tbl:synopsis of ff}
\end{table}

Before returning to the factorial and its non-recursive definition in the next section, there is an interesting application of \textit{F.F.} we can examine. It's about making an output table of three columns, like before, but the output pattern is read as one number in binary.

For example, the output pattern of $\xi^\prime(n)$ is $[1, 0, 0]$. That can be read as the binary number $100_2$, which is equal to 4 in base 10. Many \textit{F.F.} exhibit this behavior, but those with outputs like $-1$ or $n$ are excluded.

In this way, we can construct a binary table. The \ref{tbl:binary of ff} has all 3-digit binary numbers:
\begin{table}[H]
    \centering
    \begin{tabular}{c|ccc}
        $n$ & $-$ & 0 & $+$ \\
        \hline
        0 & 0 & 0 & 0 \\
        $\xi^\prime(-n)$ & 0 & 0 & 1 \\
        $\Theta(n)$ & 0 & 1 & 0 \\
        $\xi(n)$ & 0 & 1 & 1 \\
        $\xi^\prime(n)$ & 1 & 0 & 0 \\
        $Q^\prime(n)$ & 1 & 0 & 1 \\
        $\xi(-n)$ & 1 & 1 & 0 \\
        1 & 1 & 1 & 1
    \end{tabular} 
    \caption{\textit{F.F.} outputs in binary}
    \label{tbl:binary of ff}
\end{table}

The completion of this table helped define some of these functions. After discovering a few of them, there was motivation to finish \ref{tbl:binary of ff}. However, it's important to note that the binary representation doesn't have any applications for this paper; it just helps organize the table.

Since there are many new functions introduced in this paper, it's easy to get them mixed up. Refer to the previous equation list to be reminded of what each function does. That list will ensure that everything in the paper is clear and easy to understand.

This concludes the second section about \textit{foundational functions} and it's time to return to the factorial and tackle the problem of finding global, non-recursive $\prod$-product definitions of the Roman factorial.

\section{\texorpdfstring{$\prod$-product generalizations of $\rf{n}$}{Pi-product generalizations of the Roman factorial}}\label{s:rf prod gen}
\subsection{Introduction}\label{ss:rf prod gen intro}

This section is about generalizing the $\prod$-product definitions of the Roman factorial, in the same manner as it was done in \Cref{s:rf rec rl gen}.

In \Cref{ss:tf}, we introduced two formulations of the factorial function: one as a rising product and the other as a falling product. These expressions are reiterated below:
\begin{gather}
    n! = \prod_{k = 1}^{n} k \,= 1 \cdot 2 \cdot 3 \cdots (n - 2)(n - 1) \:\! n \,, \M n \in \mbZ^+ , \tag{\ref{eq:tf prod rising}} \\
    n! = \prod_{k = 0}^{n - 1} (n - k) \,= n(n - 1) \cdots 3 \cdot 2 \cdot 1 \,, \M n \in \mbZ^+ . \tag{\ref{eq:tf prod falling}}
\end{gather}

These definitions are restricted to positive integers, but we aim to extend them to encompass all $n \in \mbZ$. Our goal in this section is to derive similar factorial definitions using rising and falling products that also apply to negative integers.

We will follow a methodology akin to that presented in \Cref{s:rf rec rl gen}: utilizing \textit{F.F.}, we will formulate a piece-wise definition of the factorial. However, at present, we lack a $\prod$-product representation for negative Roman factorials, either as an ascending or descending product.

Our initial focus will be on the rising product, as depicted in \ref{eq:tf prod rising}. We aim to devise a comparable expression for negative integers while also addressing the $n = 0$ case, thus crafting a piece-wise definition. Through incremental steps, we will unify this definition using \textit{F.F.}, arriving at a comprehensive expression.

These procedures will be repeated for the falling product. The outcomes will then be summarized, therefore concluding \href{https://arxiv.org/abs/2403.09581}{Part 1} of this study.

\subsection{Rising product expression}\label{ss:rf prod neg rising}

Given the absence of a product definition for negative integers of the Roman factorial, let's initiate the construction of one by making several observations. Firstly, let's recall the values of negative Roman factorials:
\begin{table}[H]
    \centering
    \renewcommand{\thetable}{\ref{tbl:rf values}}
    \begin{tabular}{c|ccccccc}
        $n$ & -7 & -6 & -5 & -4 & -3 & -2 & -1 \\
        \hline
        $\rf{n}$ & $\nicefrac{1}{720}$ & -$\nicefrac{1}{120}$ & $\nicefrac{1}{24}$ & -$\nicefrac{1}{6}$ & $\nicefrac{1}{2}$ & -1 & 1
    \end{tabular}
    \caption{Negative Roman factorials}
    \renewcommand{\thetable}{\text{Tbl.} \arabic{section}.\arabic{subsection}.\arabic{table}}
\end{table}

It's notable that negative Roman factorials are fractions with alternating signs and a traditional factorial in the denominator. These outcomes can be reproduced by multiplying inverse negative integers, as demonstrated below:
\begin{gather*}
    \rf{-3} = \cfrac{1}{2} = \cfrac{1}{-1} \cdot \cfrac{1}{-2} \;, \\[0.618034em]
    \rf{-4} = -\cfrac{1}{6} = \cfrac{1}{-1} \cdot \cfrac{1}{-2} \cdot \cfrac{1}{-3} \;, \\[0.618034em]
    \rf{-5} = \cfrac{1}{24} = \cfrac{1}{-1} \cdot \cfrac{1}{-2} \cdot \cfrac{1}{-3} \cdot \cfrac{1}{-4} \;.
\end{gather*}

Indeed, negative Roman factorials seem to be the product of consecutive reciprocals of integers, albeit with an offset. Despite the fractions getting smaller with each increment, since the terms start at $-1$ and end at $n$, we designate this a rising product. Disregarding the offset for now, this pattern can be expressed as follows:
\begin{equation}\label{eq:rf prod neg rising wrong}
    \rf{n} = \prod_{k = 1}^{-n} \cfrac{1}{-k} \;, \M n \in \mbZ^- .
\end{equation}

However, this expression doesn't accurately represent the desired values. It appears that there is one additional multiplication term, which we refer to as an offset. For instance, \ref{eq:rf prod neg rising wrong} yields the following factorials:
\begin{gather*}
    \rf{-3} = \cfrac{1}{-1} \cdot \cfrac{1}{-2} \cdot \cfrac{1}{-3} = -\cfrac{1}{6} \:, \\[0.618034em]
    \rf{-4} = \cfrac{1}{-1} \cdot \cfrac{1}{-2} \cdot \cfrac{1}{-3} \cdot \cfrac{1}{-4} = \cfrac{1}{24} \:, \\[0.618034em]
    \rf{-5} = \cfrac{1}{-1} \cdot \cfrac{1}{-2} \cdot \cfrac{1}{-3} \cdot \cfrac{1}{-4} \cdot \cfrac{1}{-5} = -\cfrac{1}{120} \:.
\end{gather*}

These values don't match those in \ref{tbl:rf values}. As previously noted, there's an extra multiplication term at the end, which, if canceled, makes \ref{eq:rf prod neg rising wrong} accurate.

There are two approaches to address this issue: one is to adjust the upper limit of the product, subtracting 1 from it to multiply by one less term. The other involves multiplying the existing product by $n$, which cancels out the term causing the offset. Both methods work equally well, although the first idea utilizes the empty product\footnote{For more information about the empty product, check \hyperref[ss:add product]{Addendum~\ref*{ss:add product}}.} when $n = -1$.

We will proceed with the second approach for now. Therefore, multiplying each negative Roman factorial by $n$ results in the following outcomes:
\begin{gather*}
    \rf{-3} = (-3) \cdot \cfrac{1}{-1} \cdot \cfrac{1}{-2} \cdot \cfrac{1}{-3} = \cfrac{1}{2} \:, \\[0.618034em]
    \rf{-4} = (-4) \cdot \cfrac{1}{-1} \cdot \cfrac{1}{-2} \cdot \cfrac{1}{-3} \cdot \cfrac{1}{-4} = -\cfrac{1}{6} \:, \\[0.618034em]
    \rf{-5} = (-5) \cdot \cfrac{1}{-1} \cdot \cfrac{1}{-2} \cdot \cfrac{1}{-3} \cdot \cfrac{1}{-4} \cdot \cfrac{1}{-5} = \cfrac{1}{24} \:.
\end{gather*}

Thus, we formulate the following definition for negative Roman factorials as described by a rising product:
\begin{equation}\label{eq:rf prod neg rising}
    \rf{n} = n \cdot \prod_{k = 1}^{-n} \cfrac{1}{-k} \:, \M n \in \mbZ^- .
\end{equation}

Notice the similarities between this expression and \ref{eq:tf prod rising}. We will eventually consider the case of $n = 0$ and form the piece-wise definition, although let's firstly look at the analogous expression for the falling product.

\subsection{Falling product expression}\label{ss:rf prod neg falling}

In this subsection, we will derive a definition akin to \ref{eq:rf prod neg rising}, but expressed as a falling product. To start, let's reiterate the factorials of interest, this time presented as products of reciprocals of negative integers in reverse order:
\begin{gather*}
    \rf{-3} = \cfrac{1}{2} = \cfrac{1}{-2} \cdot \cfrac{1}{-1} \:, \\[0.618034em]
    \rf{-4} = -\cfrac{1}{6} = \cfrac{1}{-3} \cdot \cfrac{1}{-2} \cdot \cfrac{1}{-1} \:, \\[0.618034em]
    \rf{-5} = \cfrac{1}{24} = \cfrac{1}{-4} \cdot \cfrac{1}{-3} \cdot \cfrac{1}{-2} \cdot \cfrac{1}{-1} \:.
\end{gather*}

From our previous analysis, we observed that the negative Roman factorial of $n$ can be represented by a product consisting of $(n - 1)$ terms. Thus, the upper and lower limits of the product must be set accordingly. A possible definition expressing the aforementioned products is as follows:
\begin{equation}\label{eq:rf prod neg falling wrong}
    \prod_{k = 1}^{-n - 1} \cfrac{1}{n + k} = \cfrac{1}{n + 1} \cdot \cfrac{1}{n + 2} \cdots \cfrac{1}{-2} \cdot \cfrac{1}{-1} \:, \M n \in \mbZ^- .
\end{equation}

While accurate and applicable to all negative integers, for $n = -1$, the use of the empty product is necessary again as the upper limit becomes smaller than the lower one. It is preferable to have a definition with the same lower product limit $k = 0$ as in \ref{eq:tf prod falling}. So, we can multiply by one extra term and subsequently cancel it out, as demonstrated below:
\begin{equation}\label{eq:rf prod neg falling example}
    n \cdot \prod_{k = 0}^{-n - 1} \cfrac{1}{n + k} = \cfrac{n}{n} \cdot \cfrac{1}{n + 1} \cdots \cfrac{1}{-2} \cdot \cfrac{1}{-1} \:, \M n \in \mbZ^- .
\end{equation}

\ref{eq:rf prod neg falling example} produces the desired values while also avoiding using the empty product when $n = -1$. Therefore, we will adopt this definition for the falling product expression:
\begin{equation}\label{eq:rf prod neg falling}
    \rf{n} = n \cdot \prod_{k = 0}^{-n - 1} \cfrac{1}{n + k} \;, \M n \in \mbZ^- .
\end{equation}

In summary, in \Cref{ss:rf prod neg rising,ss:rf prod neg falling}, we've established two expressions describing negative Roman factorials as rising or falling products. They are consolidated below:
\begin{equation}\label{eq:rf prod neg rising and falling}
    \rf{n} = 
    \begin{cases}
        n \cdot \displaystyle\prod_{k = 1}^{-n} \cfrac{1}{-k} \\[1.618034em]
        n \cdot \displaystyle\prod_{k = 0}^{-n - 1} \cfrac{1}{n + k}
    \end{cases}
     , \| n \in \mbZ^- .
\end{equation}

In the subsequent subsections, we will unify the relationships for both positive and negative integers, in both the rising and falling products.

\subsection{\texorpdfstring{$\prod$-product definition step 1: $|n|$}{Pi-product definition step 1: Absolute value of n}}\label{ss:rf prod gen step 1}

Let's summarize what we have found so far about Roman factorials for positive and negative integers. Expressed as a rising product, we have:
\begin{equation}\label{eq:rf prod gen rising cases}
    \rf{n} = 
    \begin{cases}
        \,\displaystyle\prod_{k = 1}^{n} k & , \| n \in \mbZ^+ \\[1.618034em]
        n \cdot \displaystyle\prod_{k = 1}^{-n} \cfrac{1}{-k} & , \| n \in \mbZ^- .
    \end{cases}
\end{equation}

And as a falling product, we found:
\begin{equation}\label{eq:rf prod gen falling cases}
    \rf{n} = 
    \begin{cases}
        \,\displaystyle\prod_{k = 0}^{n - 1} (n - k) & , \| n \in \mbZ^+ \\[1.618034em]
        n \cdot \displaystyle\prod_{k = 0}^{-n - 1} \cfrac{1}{n + k} & , \| n \in \mbZ^- .
    \end{cases}
\end{equation}

Although the above definitions are concise, presenting both of them simultaneously can provide a clearer understanding of the generalization process. So, we present a table representing both cases of the rising and falling product expressions:
\begin{table}[h]
    \centering
    \begin{tabular}{c|cc}
        $\rf{n}$ & Rising product & Falling product \\
        \hline \\[-1.618034em]
        $n \in \mbZ^+$ & $\,\displaystyle\prod_{k = 1}^{n} k$ & $\|\displaystyle\prod_{k = 0}^{n - 1} (n - k)$ \\[1.368034em]
        $n \in \mbZ^-$ & $n \cdot \displaystyle\prod_{k = 1}^{-n} \cfrac{1}{-k}$ & $n \cdot \displaystyle\prod_{k = 0}^{-n - 1} \cfrac{1}{n + k}$
    \end{tabular}
    \caption{Roman factorial as a rising or\\falling product (not including $n = 0$)}
    \label{tbl:rf prod gen rising and falling cases without 0}
\end{table}

We will update this table as we add \textit{F.F.} into it. Our goal is to reach two \textit{universal definitions of the Roman factorial}: as an ascending and as a descending product. \ref{tbl:rf prod gen rising and falling cases without 0} serves as a way of visualising the generalization, as well as being concise.

Now, let's consider the case when $n = 0$. We see that the definitions for $n > 0$, if applied for $n = 0$, output 1. That is because the upper limits become smaller than the lower limits, thus resulting in the empty product, which is always defined to be 1. In this case specifically, we will keep the empty product in the definition for simplicity.

Since our definitions correctly describe the case of $n = 0$, we can simply extend the domain of the $n \in \mbZ^+$ cases to $n \in \mbZ^+_0$. Therefore, the previous table is updated below:
\begin{table}[H]
    \centering
    \begin{tabular}{c|cc}
        $\rf{n}$ & Rising product & Falling product \\
        \hline \\[-1.618034em]
        $n \in \mbZ^+_0$ & $\,\displaystyle\prod_{k = 1}^{n} k$ & $\|\displaystyle\prod_{k = 0}^{n - 1} (n - k)$ \\[1.368034em]
        $n \in \mbZ^-$ & $n \cdot \displaystyle\prod_{k = 1}^{-n} \cfrac{1}{-k}$ & $n \cdot \displaystyle\prod_{k = 0}^{-n - 1} \cfrac{1}{n + k}$
    \end{tabular}
    \caption{Roman factorial as a\\rising or falling product}
    \label{tbl:rf prod gen rising and falling cases}
\end{table}

Next, we notice the upper limits of the products. For $n \in \mbZ^+_0$, they are either $n$ or $n - 1$, and for $n \in \mbZ^-$ they are $-n$ or $-n - 1$.

Essentially, when $n$ becomes negative, a negative sign is applied to ensure that the upper limits remain positive numbers. As seen in \Cref{ss:rf gen step 3}, we can replace occurrences of $n$ and $-n$ with $|n|$. Consequently, the table is modified as follows:
\begin{table}[H]
    \centering
    \begin{tabular}{c|cc}
        $\rf{n}$ & Rising product & Falling product \\
        \hline \\[-1.618034em]
        $n \in \mbZ^+_0$ & $\,\displaystyle\prod_{k = 1}^{|n|} k$ & $\|\displaystyle\prod_{k = 0}^{|n|- 1} (n - k)$ \\[1.368034em]
        $n \in \mbZ^-$ & $n \cdot \displaystyle\prod_{k = 1}^{|n|} \cfrac{1}{-k}$ & $n \cdot \displaystyle\prod_{k = 0}^{|n| - 1} \cfrac{1}{n + k}$
    \end{tabular}
    \caption{Roman factorial as a rising or\\falling product (generalization: step 1)}
    \label{tbl:rf prod gen rising and falling cases step 1}
\end{table}

Unlike \Cref{s:rf rec rl gen}, this section generalizes two definitions simultaneously. While it's possible to separate the procedure for the rising product from the falling one, the steps are identical: combining the two Roman factorial product definitions into a single table allows us to present the generalization process concurrently, unifying the expressions top to bottom.

\subsection{\texorpdfstring{$\prod$-product definition step 2: $\Phi(n)$}{Pi-product definition step 2: Phi of n}}\label{ss:rf prod gen step 2}

In this step, we will introduce the function $\Phi(n)$. As defined in \Cref{ss:Phi function}, it takes the form:
\begin{equation}\tag{\ref{eq:Phi}}
    \Phi(n) = \bp{n + \Theta(n)}^{\xi^\prime(n)} =
    \begin{cases}
        1 \,, & n \in \mbZ^+_0 \\
        n \,, & n \in \mbZ^- .
    \end{cases}
\end{equation}

Note that when $n$ is negative in \ref{tbl:rf prod gen rising and falling cases step 1}, the products are multiplied by a factor of $n$, which is absent in the products describing the $n \in \mbZ^+_0$ case. This resembles the output pattern of $\Phi(n)$, which is $[n, 1, 1]$. Therefore, we can update \ref{tbl:rf prod gen rising and falling cases step 1} as such:
\begin{table}[H]
    \centering
    \begin{tabular}{c|cc}
        $\rf{n}$ & Rising product & Falling product \\
        \hline \\[-1.618034em]
        $n \in \mbZ^+_0$ & $\Phi(n) \cdot \displaystyle\prod_{k = 1}^{|n|} k$ & $\Phi(n) \cdot \displaystyle\prod_{k = 0}^{|n|- 1} (n - k)$ \\[1.368034em]
        $n \in \mbZ^-$ & $\Phi(n) \cdot \displaystyle\prod_{k = 1}^{|n|} \cfrac{1}{-k}$ & $\Phi(n) \cdot \displaystyle\prod_{k = 0}^{|n| - 1} \cfrac{1}{n + k}$
    \end{tabular}
    \caption{Roman factorial as a rising or\\falling product (generalization: step 2)}
    \label{tbl:rf prod gen rising and falling cases step 2}
\end{table}

\subsection{\texorpdfstring{$\prod$-product definition step 3: $\theta(n)$}{Pi-product definition step 3: theta of n}}\label{ss:rf prod gen step 3}

Lastly, we will address the index terms of the products, which are the expressions containing the variable $k$. We list the index terms in the next table:
\begin{table}[H]
    \centering
    \begin{tabular}{c|cc}
        Index term & Rising product & Falling product \\
        \hline
        $n \in \mbZ^+_0$ & $k$ & $n - k$ \\
        $n \in \mbZ^-$ & $\cfrac{1}{-k}$ & $\cfrac{1}{n + k}$
    \end{tabular}
    \caption{Index terms of the products in \ref{tbl:rf prod gen rising and falling cases step 2}}
    \label{tbl:rf prod gen rising and falling cases index terms}
\end{table}

Let's highlight the sign of $k$ and the exponent of the whole expression in all cases, so that we can find to find an appropriate \textit{F.F.} that describes them. For the falling product, we will keep a negative sign, like so:
\begin{table}[H]
    \centering
    \begin{tabular}{c|cc}
        Index term & Rising product & Falling product \\
        \hline
        $n \in \mbZ^+_0$ & $(+1)k^{(+1)}$ & $\bp{n - (+1)k}^{(+1)}$ \\
        $n \in \mbZ^-$ & $(-1)k^{(-1)}$ & $\bp{n - (-1)k}^{(-1)}$
    \end{tabular}
    \caption{Index terms of the products in\\\ref{tbl:rf prod gen rising and falling cases step 2}, highlighting signs and exponents}
    \label{tbl:rf prod gen rising and falling cases index terms highlighted}
\end{table}

It is evident that the sign of $k$, as well as the exponent, is represented by the same function. That function has the output pattern $[-1, 1, 1]$ and a \textit{F.F.} that matches this requirement is $\theta(n)$:
\begin{equation}\tag{\ref{eq:theta}}
    \theta(n) = \cfrac{\delta(n)}{|\delta(n)|} =
    \begin{cases}
        +1 \,, \M n \in \mbZ^+_0 \\
        -1 \,, \M n \in \mbZ^- .
    \end{cases}
\end{equation}

Thus, if we incorporate $\theta(n)$ into the previous table, we achieve a unification:
\begin{table}[H]
    \centering
    \begin{tabular}{c|cc}
        Index term & Rising product & Falling product \\
        \hline
        $n \in \mbZ$ & $\theta(n) \cdot k^{\theta(n)}$ & $\bp{n - \theta(n) \cdot k}^{\theta(n)}$
    \end{tabular}
    \caption{Index terms of the\\products in \ref{tbl:rf prod gen rising and falling cases step 2}, unified}
    \label{tbl:rf prod gen rising and falling cases index terms with ff}
\end{table}

This modification fulfills our objective of unifying the Roman factorial definitions. One last step before we proceed: we can implement a few adjustments to improve readability for the expressions regarding the index terms. They are displayed below:
\begin{gather*}
    \theta(n) \cdot k^{\theta(n)} \rightarrow \bp{k \,\theta(n)}^{\theta(n)} \\
    \bp{n - \theta(n) \cdot k}^{\theta(n)} \rightarrow \bp{n - k \,\theta(n)}^{\theta(n)}
\end{gather*}

Finally, we present the Roman factorial universally defined across all integers, described as a rising or falling product:
\begin{table}[H]
    \centering
    \begin{tabular}{c|cc}
        $\rf{n}$ & $n \in \mbZ$ \\
        \hline \\[-1.618034em]
        Rising product & $\Phi(n) \cdot \displaystyle\prod_{k = 1}^{|n|} \bp{k \,\theta(n)}^{\theta(n)}$ \\[1.368034em]
        Falling product & $\Phi(n) \cdot \displaystyle\prod_{k = 0}^{|n| - 1} \bp{n - k \,\theta(n)}^{\theta(n)}$
    \end{tabular}
    \caption{Roman factorial as a rising\\or falling product (generalized)}
    \label{tbl:rf prod gen rising and falling}
\end{table}

This result concludes \Cref{s:rf prod gen}, as well as the first part of this study.

\ref{tbl:rf prod gen rising and falling} is not unique, as there are many equivalent definitions that express the Roman factorial. An alternative formulation is provided below:
\begin{table}[H]
    \centering
    \begin{tabular}{c|cc}
        $\rf{n}$ & $n \in \mbZ$ \\
        \hline \\[-1.618034em]
        Rising product & $\displaystyle\prod_{k = 1}^{|n| - \xi^\prime(n)} \bp{k \,\theta(n)}^{\theta(n)}$ \\[1.368034em]
        Falling product & $\displaystyle\prod_{k = 0}^{|n| - 1 - \xi^\prime(n)} \bp{n - k \,\theta(n)}^{\theta(n)}$
    \end{tabular}
    \caption{Roman factorial as a rising or falling product (alternative generalization)}
    \label{tbl:rf prod gen rising and falling alt}
\end{table}

These alternative $\prod$-products would arise if we chose a different approach in \Cref{ss:rf prod neg rising,ss:rf prod neg falling}, when we chose to multiply the existing products by $n$ instead of altering their upper limits.

Here we see that there is a \textit{F.F.} in the limits of the $\prod$-products, a result we avoided for two reasons: simplicity of the expression and easier generalization in later parts of this study. The last reason was not explained here because of lacking context, but it will become apparent in the next part. There, the same process of generalization will apply to the multifactorial and its expansions into negative integers.

\section{Conclusions}\label{s:conclusions}
\subsection{Summary}\label{ss:summary}

In summary:
\begin{itemize}
    \item \textbf{\Cref{s:tf intro}}: We analyzed the factorial and presented the Roman factorial.
    \item \textbf{\Cref{s:basic ff}}: We introduced a set of 5 \textit{foundational functions} (\textit{F.F.}).
    \item \textbf{\Cref{s:rf rec rl gen}}: We employed the \textit{F.F.} shown in the previous section to rewrite the Roman factorial definition concisely, in its original formulation as well as its recursive form.
    \item \textbf{\Cref{s:advanced ff}}: We introduced another set of 5 \textit{F.F.}
    \item \textbf{\Cref{s:rf prod gen}}: We found non-recursive $\prod$-product definitions to describe the values outputted by the Roman factorial. The expressions are either a rising or a falling product, split into two cases each (positive and negative integers). Additionally, we used all \textit{F.F.} to unify the cases of these definitions, achieving universality.
\end{itemize}

In \Cref{s:basic ff,s:advanced ff} we defined the following \textit{foundational functions}:
\begin{align}
    \delta(n) & = \ff{n} + 0.5 & [-, +, +] \tag{\ref{eq:delta}} \\[0.368034em]
    \theta(n) & = \cfrac{\delta(n)}{\left|\delta(n)\right|} & [-1, 1, 1] \tag{\ref{eq:theta}} \\[0.618034em]
    \xi(n) & = \cfrac{1 + \theta(n)}{2} & [0, 1, 1] \tag{\ref{eq:xi}} \\[0.618034em]
    \xi^\prime(n) & = \cfrac{1 - \theta(n)}{2} & [1, 0, 0] \tag{\ref{eq:xi prime}} \\[0.618034em]
    \eta(n) & = \theta(n)^{-\cf{n} - 1} & [\spm 1, 1, 1] \tag{\ref{eq:eta}} \\[0.368034em]
    \Theta(n) & = \xi(n) \cdot \xi(-n) & [0, 1, 0] \tag{\ref{eq:Theta}} \\[0.368034em]
    Q(n) & = \theta(n) - \Theta(n) & [-1, 0, 1] \tag{\ref{eq:Q}} \\[0.368034em]
    Q^\prime(n) & = 1 - \Theta(n) & [1, 0, 1] \tag{\ref{eq:Q prime}} \\[0.368034em]
    \Psi(n) & = n + \Theta(n) & [n, 1, n] \tag{\ref{eq:Psi}} \\[0.368034em]
    \Phi(n) & = {\Psi(n)}^{\xi^\prime(n)} & [n, 1, 1] \tag{\ref{eq:Phi}}
\end{align}

\subsection{Results}\label{ss:results}

In \Cref{s:rf rec rl gen} we condensed the Roman factorial definition. It is defined as follows:
\begin{equation}\tag{\ref{eq:rf pw}}
    \rf{n} = 
    \begin{cases}
        n! & , \| n \in \mbZ^+_0 \\
        \cfrac{(-1)^{-n - 1}}{(-n - 1)!} & , \| n \in \mbZ^- ,
    \end{cases}
\end{equation}
in which the factorial is defined recursively as
\begin{equation}\tag{\ref{eq:tf rec}}
    n! = n(n - 1)! \,, \M 0! = 1 \,, \M n \in \mbZ^+ .
\end{equation}

The generalized relationship is
\begin{equation}\tag{\ref{eq:rf gen}}
    \rf{n} = \eta(n) \cdot \bp{|n| - \xi^\prime(n)}!^{\,\theta(n)} \,, \M n \in \mbZ \,,
\end{equation}
where
\begin{equation}\tag{\ref{eq:tf rec}}
    n! = n(n - 1)! \,, \M 0! = 1 \,, \M n \in \mbZ^+ .
\end{equation}

Additionally, we rewrote the following doubly-recursive definition of the Roman factorial:
\begin{equation}\tag{\ref{eq:rf rec pw}}
    \rf{n} = 
    \begin{cases}
        n \rf{n - 1} & , \| n \in \mbZ^+ \\[0.368034em]
        \cfrac{\rf{n + 1}}{n + 1} & , \| n \in \mbZ^- \setminus \{-1\} \,,
    \end{cases}
\end{equation}
where
\begin{equation}\tag{\ref{eq:rf rec seeds}}
    \rf{0} = \rf{-1} = 1 \,.
\end{equation}

The outcome of the generalization is as follows:
\begin{equation}\tag{\ref{eq:rf rec gen}}
    \rf{n} = \bp{n + \xi^\prime(n)}^{\theta(n)} \rf{n - \theta(n)} \,, \M n \in \mbZ \setminus \{0, \|-1\} \,,
\end{equation}
where
\begin{equation}\tag{\ref{eq:rf rec seeds}}
    \rf{0} = \rf{-1} = 1 \,.
\end{equation}

Lastly, in \Cref{s:rf prod gen}, we constructed these definitions:
\begin{table}[H]
    \centering
    \renewcommand{\thetable}{\ref{tbl:rf prod gen rising and falling cases}}
    \begin{tabular}{c|cc}
        $\rf{n}$ & Rising product & Falling product \\
        \hline \\[-1.618034em]
        $n \in \mbZ^+_0$ & $\,\displaystyle\prod_{k = 1}^{n} k$ & $\|\displaystyle\prod_{k = 0}^{n - 1} (n - k)$ \\[1.368034em]
        $n \in \mbZ^-$ & $n \cdot \displaystyle\prod_{k = 1}^{-n} \cfrac{1}{-k}$ & $n \cdot \displaystyle\prod_{k = 0}^{-n - 1} \cfrac{1}{n + k}$
    \end{tabular}
    \caption{Roman factorial as a\\rising or falling product}
    \renewcommand{\thetable}{\text{Tbl.} \arabic{section}.\arabic{subsection}.\arabic{table}}
\end{table}

The two definitions in the above table were consolidated into the formulations listed below:
\begin{table}[H]
    \centering
    \renewcommand{\thetable}{\ref{tbl:rf prod gen rising and falling}}
    \begin{tabular}{c|cc}
        $\rf{n}$ & $n \in \mbZ$ \\
        \hline \\[-1.618034em]
        Rising product & $\Phi(n) \cdot \displaystyle\prod_{k = 1}^{|n|} \bp{k \,\theta(n)}^{\theta(n)}$ \\[1.368034em]
        Falling product & $\Phi(n) \cdot \displaystyle\prod_{k = 0}^{|n| - 1} \bp{n - k \,\theta(n)}^{\theta(n)}$ \\[1.368034em]
    \end{tabular}
    \caption{Roman factorial as a rising or\\falling product (generalized)}
    \renewcommand{\thetable}{\text{Tbl.} \arabic{section}.\arabic{subsection}.\arabic{table}}
\end{table}

In the next part, we will examine the double factorial and follow a similar procedure. We will investigate an extension into negative integers, we will find recursive as well as non-recursive piece-wise definitions of the expanded double factorial and we will unify their cases. This process will be repeated for the triple factorial and in general for all factorials of higher orders.

\section{Details and References}\label{s:references}
\phantomsection

\begin{flushleft}
    \large\textbf{References}
    \addcontentsline{toc}{subsection}{References}
    \printbibliography[heading=none]

@book{roman,
    author    = {Daniel E. Loeb and Gian-Carlo Rota},
    title     = {Formal Power Series of Logarithmic Type},
    series    = {Advances in Mathematics},
    year      = {1989},
    publisher = {Elsevier},
    doi       = {10.1016/0001-8708(89)90079-0}
}

@article{0tothe0,
    ISSN      = {00029890, 19300972},
    URL       = {http://www.jstor.org/stable/2307224},
    author    = {L. J. Paige},
    journal   = {The American Mathematical Monthly},
    number    = {3},
    pages     = {189--190},
    publisher = {Mathematical Association of America},
    title     = {A Note on Indeterminate Forms},
    volume    = {61},
    year      = {1954},
    doi       = {10.2307/2307224}
}

@book{fallrisefact,
    author    = {Steffensen, Johan Frederik},
    title     = {Interpolation},
    edition   = {2nd},
    publisher = {Dover Publications},
    year      = {2006},
    month     = {3},
    isbn      = {0-486-45009-0},
    note      = {A reprint of the 1950 edition by Chelsea Publishing.}
}

@online{boolean_algebra,
    author  = {Wikipedia},
    title   = {Boolean Algebra},
    urldate = {2024-03-03},
    url     = {https://en.wikipedia.org/wiki/Boolean_algebra}
}

@online{boolean_functions,
    author  = {Wikipedia},
    title   = {Boolean-valued function},
    urldate = {2024-03-03},
    url     = {https://en.wikipedia.org/wiki/Boolean-valued_function}
}

@online{0_factorial,
    author  = {Wikipedia},
    title   = {Factorial},
    url     = {https://en.wikipedia.org/wiki/Factorial},
    urldate = {2024-03-07}
}

@misc{romanfact,
    author        = {Daniel E. Loeb},
    title         = {A generalization of the binomial coefficients}, 
    year          = {1995},
    eprint        = {math/9502218},
    archivePrefix = {arXiv},
    primaryClass  = {math.CO}, 
    doi           = {10.48550/arXiv.math/9502218}
}

@online{product,
    author  = {Wikipedia},
    title   = {Product (mathematics)},
    urldate = {2024-03-13},
    url     = {https://en.wikipedia.org/wiki/Product_(mathematics)}
}

@online{setwiki,
    author  = {Wikipedia},
    title   = {Set (mathematics)},
    urldate = {2024-03-13},
    url     = {https://en.wikipedia.org/wiki/Set_(mathematics)}
}

@online{setfun,
    author  = {Math is Fun},
    title   = {Common number sets},
    urldate = {2024-03-13},
    url     = {https://www.mathsisfun.com/sets/number-types.html},
    note    = {Blog}
}

@book{titchmarsh1986theory,
  title     = {The Theory of the Riemann Zeta-function},
  author    = {Titchmarsh, E.C. and Heath-Brown, D.R.},
  isbn      = {9780198533696},
  lccn      = {lc86012520},
  series    = {Oxford science publications},
  url       = {https://books.google.gr/books?id=1CyfApMt8JYC},
  year      = {1986},
  publisher = {Clarendon Press}
}
\end{flushleft}
\newpage

\detnref{Acknowledgements}
{Special thanks go to Konstantinos Tsakalidis, my highschool math teacher, who realized my potential, engaged in my curiosity and urged me to continue my research.}

\detnref{Author's Contributions}
{The author has conceived the ideas, made the calculations and has also written and approved the manuscript. The first version of this work was submitted to \href{https://arxiv.org/abs/2403.09581}{\ul{arXiv.org}} on $\pi$ day: March 14, 2024.}

\detnref{Software Used}
{This document was made in \LaTeX \,using the online platform Overleaf, compiled in \TeX \text{ }Live 2023. All figures were drawn in Mathematica. ChatGPT 3.5 has assisted in the formatting of this document and \ul{not} in any research presented here.}

\detnref{Code Availability}
{The Mathematica code used to generate the figures and perform the calculations in this document is available upon request. The \LaTeX \, document is available for downloading on \href{https://arxiv.org/abs/2403.09581}{\ul{arXiv.org}}.

For any questions, feedback, error correction or further discussion, feel free to contact me at this email: \href{mailto:lliponis@physics.auth.gr}{\ul{lliponis@physics.auth.gr}}.}

\section{Addendum}\label{s:addendum}
\subsection{Number sets}\label{ss:add number sets}

A set is a mathematical concept representing a collection of distinct items, called elements or members \cite{setwiki} \cite{setfun}. Some sets are so significant in mathematics that they have special names and notations to identify them. These important sets are often denoted using blackboard bold typeface (e.g. $\mbZ$). Common number sets include the following sets:
\begin{table}[H]\label{tbl:add number sets}
    \centering
    \begin{tabular}{c|cc}
        Set & Description & Example \\
        \hline
        $\mbN \,, \mbZ^+$ & Natural numbers & 1, 2, 3, 4 $\cdots$ \\
        $\mbN_0$ & $\mbN$ with 0 (or $\mbZ^+_0)$ & 0, 1, 2, 3 $\cdots$ \\
        $\mbZ$ & Integers & $\cdots$ -2, -1, 0, 1 $\cdots$ \\
        $\mathbb{Q}$ & Rational numbers & \nicefrac{1}{2}, -\nicefrac{5}{4}, 0.01 $\cdots$ \\
        $\mbR$ & Real numbers & $\sqrt{2}$, $\pi$, $e$, $\phi$ $\cdots$ \\
        $\mathbb{I}$ & Imaginary numbers & $i$, $9.7i$, $-\nicefrac{i}{25}$ $\cdots$ \\
        $\mathbb{C}$ & Complex numbers & $1 + i$, $\sqrt{3} - 6i$ $\cdots$
    \end{tabular}
    \caption{Common number sets}
\end{table}

Note that $\mathbb{C}$ includes all real and imaginary numbers. For instance, the number $2 - 3i$ is itself a complex number, as well as $5.2i$ and $4$.

By combining these number sets, we can express other sets. For instance, we can represent the set of all real numbers excluding negative integers as
\begin{equation*}
    \mbR \setminus \mathbb{Z^-} \,,
\end{equation*}
where $\mathbb{Z^-}$ represents negative integers.

Also, number sets with an asterisk (*) usually denote that they do not include the number 0. For instance, $\mbZ^*$ refers to all integers except 0.

Additionally, a number set without a collection of elements can be expressed in detail. For example, the set of all positive odd integers except 1, is:
\begin{equation*}
    \mbZ^+_{odd} \setminus \{1\} \,.
\end{equation*}

\subsection{\texorpdfstring{$\prod$-product}{Pi-product}}\label{ss:add product}

The product operator $\prod$ for a product of a sequence is represented by the capital Greek letter "pi" ($\Pi$), similar to how $\sum$ is used as the summation symbol \cite{product}. For example, the product of the first 6 squares of natural numbers can be expressed as:
\begin{equation*}
    \prod_{k = 1}^{6} k^{2} = 1 \cdot 4 \cdot 9 \cdot 16 \cdot 25 \cdot 36 \,.
\end{equation*}

The number above the symbol $\prod$ is called the upper limit of the product, while the number below is the lower limit. In the example provided, the upper limit is 6, and the lower limit is 1. The variable $k$ denotes the multiplicands or factors of the product.

Moreover, if the terms of the product increase successively, it's called a rising product. Conversely, if they decrease, it's referred to as a falling product.

For the $\prod$-product to be well-defined, it's typically required that both the upper and lower limits be integers, often natural numbers.

When both limits are set to a specific number, the product is evaluated to that particular number:
\begin{equation*}
    \displaystyle\prod_{k = 5}^5 k = 5 \,, \|\prod_{k = 3}^3 e^k = e^3.
\end{equation*}

If there are no factors at all, it results in what's known as the empty product, which is defined as 1. This occurs when the upper limit of the product is less than the lower limit by any amount, regardless of the index expression of $k$. For instance:
\begin{equation*}
    \displaystyle\prod_{k = 3}^2 k = \prod_{k = 3}^2 2k = \prod_{k = 3}^2 k^3 = 1 \,.
\end{equation*}

In summary, we have the following identities:
\begin{equation}\label{eq:add prod properties}
    \displaystyle\prod_{k = n}^n f(k) = f(n) \,, \|\displaystyle\prod_{k = n}^{n - a} f(k) = 1 \,, \|a \in \mbZ^+ , \| n \in \mbZ \,,
\end{equation}
where $f(k)$ is an arbitrary function of $k$.

\subsection{Falling and rising factorials}\label{ss:add fall rise fact}

The falling factorial (also known as the descending factorial, falling sequential product, or lower factorial \cite{fallrisefact}) is defined as the polynomial
\begin{equation}\label{eq:add fall fact}
    (x)_n = \prod_{k = 0}^{n - 1}(x - k) = x(x - 1)(x - 2) \cdots (x - n + 1) \,.
\end{equation}

The rising factorial (also known as the Pochhammer function, Pochhammer polynomial, ascending factorial, rising sequential product, or upper factorial \cite{fallrisefact}) is defined as
\begin{equation}\label{eq:add rise fact}
    x_{(n)} = \prod_{k = 0}^{n - 1}(x + k) = x(x + 1)(x + 2) \cdots (x + n - 1) \,.
\end{equation}

Each of these symbols is defined to be 1 (representing an empty product) when $n = 0$. Together, they are referred to as factorial powers. They are directly related to the ordinary factorial:
\begin{gather}
    n! = 1^{(n)} = (n)_n \,, \label{eq:add rise fall fact property 1} \\
    (x)_n = \cfrac{x!}{(x - n)!} \:, \label{eq:add rise fall fact property 2} \\
    x^{(n)} = \cfrac{(x + n - 1)!}{(x - 1)!} \:. \label{eq:add rise fall fact property 3}
\end{gather}

Falling and rising factorials should not be confused with falling and rising products: the latter are any products where the terms themselves either increase or decrease. We will not discuss falling or rising factorials in this paper.

\subsection{\texorpdfstring{The function $\delta(n)_f$}{The function Fourier delta of n}}\label{ss:add ff deltaf}

The function $\delta(n)$, defined in \Cref{ss:delta function}, is expressed as:
\begin{equation}\tag{\ref{eq:delta}}
    \delta(n) = \ff{n} + \varepsilon \,, \M 0 < \varepsilon < 1 \,, \M n \in \mbR \,.
\end{equation}

The primary objective of $\delta(n)$ is to yield a positive value for non-negative inputs and a negative value for negative inputs.

An alternate representation of this function involves a Fourier approximation. However, before presenting that, it is necessary to examine the sawtooth function.

Essentially, every number can be split into two parts: the integer part (called the floor function $\ff{n}$) and the fractional part (referred to as the sawtooth function $\{n\}$):
\begin{equation}\label{eq:n decomposition as floor and sawtooth}
    n = \ff{n} + \{n\} \,, \M 0 \leq \{n\} < 1 \,.
\end{equation}

To comprehend this decomposition, consider the following examples:
\begin{gather*}
    5.7 = \ff{5.7} + \{5.7\} = 5 + 0.7 \,, \\
    6 = \ff{6} + \{6\} = 6 + 0 \,, \\
    -0.2 = \ff{-0.2} + \{-0.2\} = -1 + 0.8 \,, \\
    -3.9 = \ff{-3.9} + \{-3.9\} = -4 + 0.1 \,.
\end{gather*}

It's important to note that the floor function of a number, as defined here, is the largest integer less than or equal to that number. For negative numbers, taking the floor involves rounding down to the next smaller integer (further from 0). The fractional part always lies between 0 and 1 in this context.

A formal mathematical definition of the floor function is provided as well:
\begin{equation}\label{eq:floor def}
    \ff{x} = max\{m \in \mbZ \|| \|m \leq x \} \,.
\end{equation}

The sawtooth function $\{n\}$ is graphically represented as follows:
\begin{figure}[H]
    \centering
    \includegraphics[width = 1\linewidth]{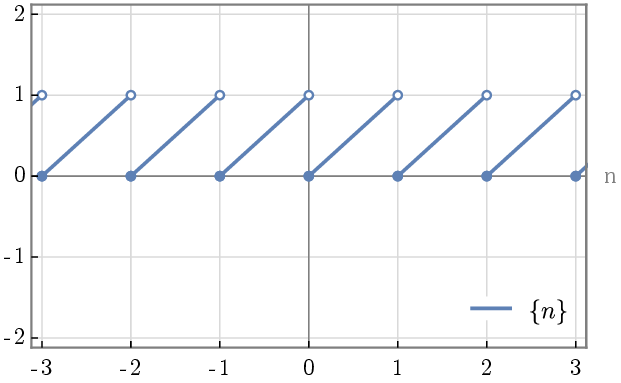}
    \caption{The sawtooth function $\{n\}$}
    \label{fig:sawtooth}
\end{figure}

Now, let's explore the Fourier approximation of the sawtooth function \cite{titchmarsh1986theory}. This approximation, denoted $\{n\}_f$ with the subscript $f$ for Fourier, is defined as:
\begin{equation}\label{eq:fourier sawtooth}
    \{n\}_f = \frac{1}{2} - \frac{1}{\pi} \sum_{k = 1}^{\infty} \frac{\sin(2\pi kn)}{k} \:, \M n \in \mbR \,.
\end{equation}

At points of discontinuity, a Fourier approximation equals the midpoint between the discontinuity points. In this case, \ref{eq:fourier sawtooth} is equal to $\nicefrac{1}{2}$ at integers: $\{0\}_f = \{1\}_f = \nicefrac{1}{2}$. Thus, the plot of $\{n\}_f$ resembles this:
\begin{figure}[H]
    \centering
    \includegraphics[width = 1\linewidth]{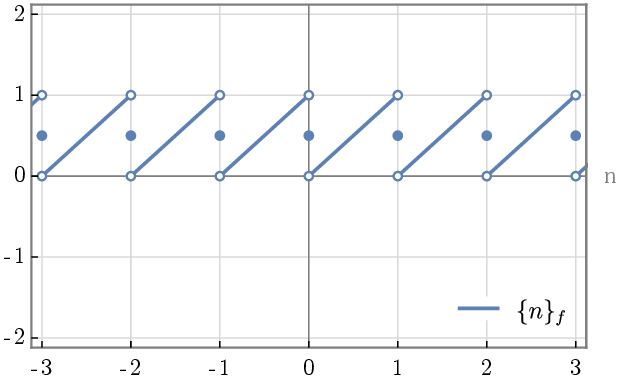}
    \caption{The sawtooth function $\{n\}_f$}
    \label{fig:fourier sawtooth}
\end{figure}

The next step is to incorporate this Fourier-approximated floor function, denoted as $\ff{n}_f$, into the definition of $\delta(n)$. This new expression is obtained by substituting $\{n\}$ for $\{n\}_f$ in \ref{eq:n decomposition as floor and sawtooth}:
\begin{equation}\label{eq:fourier floor}
    \ff{n}_f = n - \{n\}_f = n - \frac{1}{2} + \frac{1}{\pi} \sum_{k = 1}^{\infty} \frac{\sin(2\pi kn)}{k} \:, \|\: n \in \mbR \,.
\end{equation}

The graphical representation of $\ff{n}_f$ is illustrated below:
\begin{figure}[H]
    \centering
    \includegraphics[width = 1\linewidth]{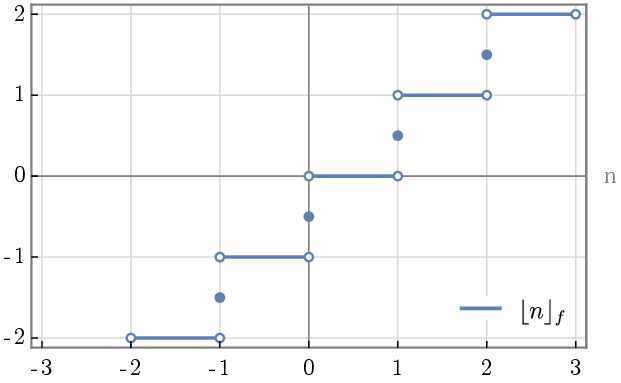}
    \caption{The floor function $\ff{n}_f$}
    \label{fig:fourier floor}
\end{figure}

We can see that it closely resembles the traditional floor function, with the distinction that the value at integers is lower by $\nicefrac{1}{2}$.

Incorporating this new floor function $\ff{n}_f$ into $\delta(n)$, we arrive at:
\begin{equation}\label{eq:delta fourier before varepsilon}
    \delta(n)_f = \ff{n}_f + \varepsilon \,, \M n \in \mbR \,,
\end{equation}
where we need to set an appropriate range for $\varepsilon$.

As previously stated in \ref{eq:delta}, $\varepsilon > 0$ ensures $\delta(0) = 0 + \varepsilon > 0$, while $\varepsilon < 1$ ensures $\delta(-1) = -1 + \varepsilon < 0$, achieving the desired output pattern $[-, +, +]$ for all real numbers.

Reconsidering this approach, we initially observe that $\delta(0)_f = -\nicefrac{1}{2}$ for $\varepsilon = 0$. To make $\delta(0)_f > 0$, $\varepsilon$ needs to be greater than $\nicefrac{1}{2}$. However, $\varepsilon$ must remain less than 1, as any value higher would result in $\delta(n) > 0$ for $n \in (-1, 0)$.

Given this analysis, we conclude that $\varepsilon \in (\nicefrac{1}{2}, 1)$. Selecting $\varepsilon = \nicefrac{3}{4}$ yields the following expression for $\delta(n)_f$:
\begin{equation}\label{eq:delta fourier}
    \delta(n)_f = n + \frac{1}{4} + \frac{1}{\pi} \sum_{k = 1}^{\infty} \frac{\sin(2\pi kn)}{k} \:, \M n \in \mbR \,.
\end{equation}

With this value of $\varepsilon$, $\delta(n)_f$ shifts positively, ensuring $\delta(-1)_f = -\nicefrac{3}{4} < 0$, $\delta(0)_f = \nicefrac{1}{4} > 0$, and $\delta(1)_f = \nicefrac{5}{4} > 0$. The graphical representation confirms this behavior:
\begin{figure}[H]
    \centering
    \includegraphics[width = 1\linewidth]{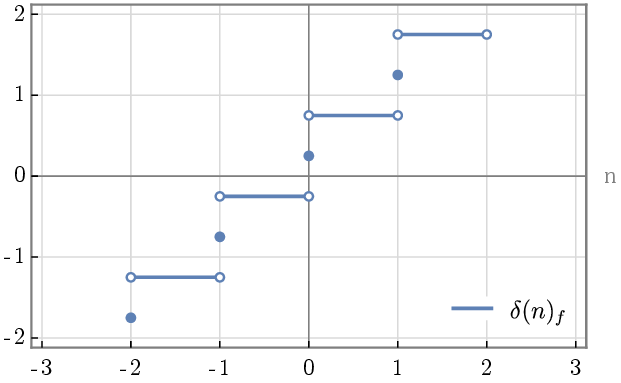}
    \caption{The function $\delta(n)_f$}
    \label{fig:delta fourier}
\end{figure}

Consequently, $\delta(n)_f$ can be used in place of $\delta(n)$ interchangeably. Either definition of \ref{eq:delta} and \ref{eq:delta fourier} will suffice, since the \textit{F.F.} $\theta(n)$ will still have the output pattern $[-1, 1, 1]$ and $\delta(n)$ does not have any other use.

\subsection{Boolean algebra}\label{ss:add boolean algebra}

Boolean algebra \cite{boolean_algebra} is a branch of algebra named after George Boole, who introduced it in his book \textit{The Mathematical Analysis of Logic} (1847). It is a formal way of describing logical operations, like logic gates which are used in computing.

In Boolean algebra, the values of the variables are the truth values \textit{true} or \textit{false}, usually denoted with 1 and 0. This field deals with functions which have their outputs in the set $\{0, 1\}$, termed Boolean-valued functions.

A Boolean-valued function \cite{boolean_functions} (sometimes called a predicate or a proposition) is a function whose input is an arbitrary number set and has outputs in a Boolean domain (i.e. a generic two-element set, for example $\{0, 1\}$). Those elements are interpreted as logical values, for example, $0 = false$ and $1 = true$ (a single bit of information).

In \Cref{s:basic ff,s:advanced ff}, we define a set of Boolean-like functions in order to unify piece-wise definitions of the Roman factorial and its various expansions and generalizations. These functions are termed here \textit{foundational functions} (\textit{F.F.}) and are similar to Boolean-valued functions.

\textit{F.F.} do not entirely belong in Boolean algebra because they have up to 3 different outputs, depending if their input $n$ is greater, less than or equal to 0. Additionally, those outputs are not constrained to be 0 or 1, but also $-1$ and even $n$ in the case of $\Phi(n)$.

However, a few \textit{F.F.} are considered to be genuine Boolean-valued functions: for instance, $\xi(n)$ is a function that outputs the true/false states of the criterion $n \geq 0$. The function $\Theta(n)$, introduced in \Cref{ss:Theta function}, can be thought of as an operation between two Boolean-valued functions:
\begin{gather*}
    (n = 0) = \big\{ (n \geq 0) \|and \|(n \leq 0) \big\} \\[1em]
    \Rightarrow\|\Theta(n) = \xi(n) \cdot \xi(-n) \,.
\end{gather*}

In this paper, there are references to Boolean algebra but further knowledge of the field is not required for understanding the premise and use of \textit{F.F.}

\subsection{\texorpdfstring{The indeterminate form $0^0$}{The indeterminate form 0 to the 0}}\label{ss:add 0 to the 0}

The mathematical expression $0^0$, known as zero to the power of zero, presents a unique case in mathematics where its definition as either 1 or left undefined depends on the context \cite{0tothe0}. In algebra and combinatorics, it is conventionally defined as $0^0 = 1$. However, in mathematical analysis, this expression is often left undefined.

Several widely used formulas involving natural-number exponents necessitate the definition of $0^0$ as 1. For instance, considering $b^0$ as an empty product yields its value as 1. Similarly, the combinatorial interpretation of $b^0$ corresponds to the count of 0-tuples from a set with $b$ elements, resulting in exactly one 0-tuple. Moreover, in set theory, $b^0$ signifies the count of functions from the empty set to a set with $b$ elements, yielding precisely one such function.

These interpretations converge to $0^0 = 1$ in various contexts. We list a few of those contexts here:

\begin{itemize}
    \item The binomial theorem $(1 + x)^n = \sum_{k = 0}^n \binom{n}{k} x^k$ holds true for $x = 0$ only if $0^0 = 1$.
    \item Similarly, in rings of power series, where $x^0$ must be defined as 1 for all specializations of $x$, identities like $\frac{1}{1 - x} = \sum_{n = 0}^\infty x^n$ and $e^x = \sum_{n = 0}^\infty \frac{x^n}{n!}$ hold true for $x = 0$ only if $0^0 = 1$.
    \item In order for the polynomial $x^0$ to define a continuous function $f: \mbR \rightarrow \mbR$, it is necessary to define $0^0 = 1$.
    \item In calculus, the power rule $\frac{d}{dx} x^n = nx^{n - 1}$ is valid for $n = 1$ at $x = 0$ only if $0^0 = 1$.
\end{itemize}

For those interested, we note a few limits of the form $0^0$ that approach different values:
\begin{gather*}
    \lim_{x \to 0^+} x^x = 1 \qquad \lim_{x \to 0^+} \bp{ \sqrt{x + 1} - \sqrt{x}}^{\nicefrac{1}{ln(lnx)}} = 0 \\[1.118034em]
    \lim_{x \to 0^+} \bp{e^{-\nicefrac{1}{x^2}}}^{-x} = +\infty \qquad \lim_{x \to 0^+} \bp{e^{-\nicefrac{1}{x}}}^{ax} = e^{-a} \\[1.118034em]
    \lim_{x \to +\infty} x^{\nicefrac{1}{ln(3x)}} = e \qquad \lim_{x \to 0^+} \bp{e^{-\nicefrac{1}{x^2}}}^x = 0
\end{gather*}

Generally, $0^0$ is defined to be 1 so that many identities and properties stay true for $n = 0$. One way of overcoming the anomaly is to use the \textit{foundational function} $\Psi(n)$ in the base, which is never 0. Thus, an expression like $\Psi(n)^n$ is well defined for $n \in \mbR$. We chose to avoid using $0^0 = 1$ in this paper since it wasn't necessary in any \textit{F.F.} or generalization.

In the examples above, the second limit in the first row is credited to a blackpenredpen's \href{https://www.youtube.com/watch?v=X65LEl7GFOw}{\ul{video}}, titled \textit{finally $0^0$ approaches 0 after 6 years!}.

\subsection{\texorpdfstring{The case of $0!$}{The case of 0 factorial}}\label{ss:add 0 factorial}

  The factorial of $0$ is defined as $1$, denoted as $0! = 1$ \cite{0_factorial}. This definition is justified by several reasons:

\begin{itemize}
    \item When $n = 0$, the definition of $n!$ involves multiplying no numbers, adhering to the convention that an empty product equals the multiplicative identity.
    \item There exists only one permutation of zero objects: since there are no objects to rearrange, the only option is to leave them as they are.
    \item Defining $0! = 1$ ensures that various combinatorial identities hold true for all valid parameter values. For example, the number of ways to select all $n$ elements from a set of $n$ is $\binom{n}{n} = \frac{n!}{n! \, 0!} = 1$, a binomial coefficient identity valid only if $0! = 1$.
    \item With $0! = 1$, the factorial's recurrence relation remains valid for $n = 1$, simplifying recursive factorials by needing only the base case of zero.
    \item Setting $0! = 1$ allows for concise expressions of many formulas, such as the exponential function represented as a power series: $e^x = \sum_{n = 0}^{\infty} \frac{x^n}{n!}\,$.
    \item This choice aligns with the Gamma function, where $0! = \Gamma(0 + 1) = 1$, a necessary condition for the gamma function to be continuous.
\end{itemize}

In this study, we do not strive to prove why/if relationships like $0! = 1$ are true, instead we try to find intuitive reasons for their validity in this context.

\end{document}